\documentclass{amsart}

\usepackage[cmtip,matrix,arrow]{xy}
\CompileMatrices
\newdir{((}{{\lhook\kern-.1em}}
\newdir{))}{{\rhook\kern-.1em}}
\newdir{ >}{{}*!/-5pt/@{>}}

\usepackage{euler,eucal}

\DeclareMathAlphabet{\mathbf}{T1}{ppl}{bx}{n}
\DeclareMathAlphabet{\mathrm}{T1}{ppl}{m}{n}

\usepackage{amssymb,amsxtra}

\usepackage{epic,eepic}

\usepackage{multicol}

\makeatletter

\def\theglossary{\@restonecoltrue\if@twocolumn\@restonecolfalse\fi
\columnseprule\z@ \columnsep 35\p@
\let\@makessectionhead\indexsec
\@xp\section\@xp*\@xp{\glossaryname}%
\let\item\@idxitem
\parindent\z@  \parskip\z@\@plus.3\p@\relax
\footnotesize}

\def\glossaryname{Notation Index}

\makeatother

\makeglossary

\numberwithin{equation}{section}
\renewcommand\dots{\relax\ifmmode\ldots\else$\,\ldots\,$\fi}

\newcommand\note[1]%
{$^\dagger$\marginpar{\footnotesize{$^\dagger${#1}}}}

\divide\count253 by 60 
\divide\count255 by 60 
\multiply\count255 by 60 
\advance\count254 by -\count255 

\def\today{\number\day\space\ifcase\month\or January\or February\or
March\or April\or May\or June\or July\or August\or September\or
October\or November\or December\fi\space\number\year}

\def\hour{\ifnum\count253<10
0\number\count253\else\number\count253\fi}

\def\minute{\ifnum\count254<10
0\number\count254\else\number\count254\fi}

\swapnumbers
\newtheorem{theorem}[equation]{Theorem}
\newtheorem{proposition}[equation]{Proposition}
\newtheorem{lemma}[equation]{Lemma}
\newtheorem{corollary}[equation]{Corollary}
\newtheorem{addendum}[equation]{Addendum}

\theoremstyle{definition}
\newtheorem{definition}[equation]{Definition}
\newtheorem{example}[equation]{Example}
\newtheorem{remark}[equation]{Remark}

\newcommand\eu{\mathfrak}
\newcommand\lie{\mathfrak}
\renewcommand\k{\lie{k}} 
\renewcommand\t{\lie{t}}
\newcommand\g{\lie{g}}
 
\newcommand\n{\lie{n}} 
\renewcommand\u{\lie{u}}

\newcommand\z{\lie{z}}

\newcommand\bb[1]{{\text{\bf#1}}}

\newcommand\Z{\bb{Z}} 

\newcommand\R{\bb{R}} 
\newcommand\C{\bb{C}}

\renewcommand\H{\bb{H}}
\renewcommand\P{\bb{P}}

\newcommand\ca{\mathcal}

\newcommand\func[1]{\operatorname{\mathrm{#1}}}
\newcommand\funclim[1]{\operatorname*{\mathrm{#1}}}

\newcommand\Ad{\func{Ad}}
\newcommand\ad{\func{ad}}

\renewcommand\cosh{\func{cosh}}

\newcommand\diag{\func{diag}}

\renewcommand\dim{\func{dim}}

\renewcommand\exp{\func{exp}}
\renewcommand\gcd{\func{gcd}}

\newcommand\Hom{\func{Hom}}
\newcommand\id{\func{id}}

\renewcommand\ker{\func{ker}}
\newcommand\lcm{\func{lcm}}
\renewcommand\lim{\funclim{lim}}
\renewcommand\log{\func{log}}

\renewcommand\min{\func{min}}
\renewcommand\mod{\func{mod}}

\newcommand\rank{\func{rank}}

\renewcommand\sin{\func{sin}}
\renewcommand\sinh{\func{sinh}}

\renewcommand\star{\func{star}}

\newcommand\tr{\func{tr}}

\renewcommand\Im{\operatorname{\mathrm{Im}}}

\newcommand\Lie{\operatorname{\mathrm{Lie}}} 

\newcommand\group[1]{{\text{\bf#1}}}

\newcommand\SO{\group{SO}}

\newcommand\SU{\group{SU}}
\newcommand\U{\group{U}}
\newcommand\Sp{\group{Sp}}

\newcommand\A{\group{A}}
\newcommand\B{\group{B}}

\newcommand\D{\group{D}}
\newcommand\E{\group{E}}
\newcommand\F{\group{F}}
\newcommand\G{\group{G}}

\newcommand\abs[1]{\lvert#1\rvert}

\newcommand\norm[1]{\lVert#1\rVert}

\newcommand\inner[1]{\langle#1\rangle}
\newcommand\biginner[1]{\bigl\langle#1\bigr\rangle}

\newcommand\quot[1][\kern.3ex]{/\kern-.7ex/_{\kern-.4ex#1}}
\newcommand\bigquot[1][\,\,]{\big/\kern-.85ex\big/_{\!\!#1}}

\newcommand\powl{[\kern-.3ex[}
\newcommand\powr{]\kern-.3ex]}
\newcommand\bigpowl{\bigl[\kern-.6ex\bigl[}
\newcommand\bigpowr{\bigr]\kern-.6ex\bigr]}

\newcommand\sur{\mathrel{\to\kern-1.8ex\to}}
\newcommand\iso{\mathrel{\hookrightarrow\kern-1.8ex\to}}

\newcommand\longhookrightarrow{\lhook\joinrel\longrightarrow}

\newcommand\longsur{\mathrel{\longrightarrow\kern-1.8ex\to}}
\newcommand\longiso{\mathrel{\longhookrightarrow\kern-1.8ex\to}}

\newcommand\alcove{\ca{A}}

\newcommand\cross{\Phi\inv(\exp\bar{\ca{A}})}
\newcommand\crosschamber{\Phi\inv(\bar{\ca{C}})}
\newcommand\dirac{/\kern-1.2ex\partial} 
\newcommand\eps{\varepsilon}
\newcommand\fusion{\circledast}

\newcommand\antiddots{\mathinner{%
\mkern1mu\raise1pt\vbox{\kern7pt\hbox{.}}
\mkern2mu\raise4pt\hbox{.}\mkern2mu\raise7pt\hbox{.}\mkern1mu}}

\newcommand\inv{^{-1}} 

\renewcommand\subset{\subseteq}

\newcommand\affine{_{\mathrm{affine}}}

\newcommand\impl{_{\mathrm{impl}}}
\newcommand\prin{_{\mathrm{prin}}}

\begin{document} 


\title{Group-valued Implosion and Parabolic Structures}

\author{Jacques Hurtubise}

\email{hurtubis@math.mcgill.ca}

\address{Department of Mathematics, McGill University, Montr\'eal, QC
H3A 2K6, and Centre de Recherches Math\'ematiques, Universit\'e de
Montr\'eal, Montr\'eal, QC H3P 3J7, Canada}

\author{Lisa Jeffrey}

\email{jeffrey@math.toronto.edu}

\address{Department of Mathematics, University of Toronto, Toronto,
Ontario M5S 3G3, Ca\-na\-da}

\author{Reyer Sjamaar}

\email{sjamaar@math.cornell.edu}

\address{Department of Mathematics, Cornell University, Ithaca, New
York 14853-7901, USA}

\thanks{J. Hurtubise was partially supported by grants from NSERC and
FQRNT, L. Jeffrey was partially supported by grants from NSERC, and
R. Sjamaar was partially supported by NSF Grant DMS-0071625.  LJ and
RS gratefully acknowledge the hospitality of the Erwin Schr\"odinger
International Institute for Mathematical Physics}

\date{27 February 2004}

\begin{abstract}
The purpose of this paper is twofold.  First we extend the notion of
symplectic implosion to the category of quasi-Hamiltonian
$K$-manifolds, where $K$ is a simply connected compact Lie group.  The
imploded cross-section of the double $K\times K$ turns out to be
universal in a suitable sense.  It is a singular space, but some of
its strata have a nonsingular closure.  This observation leads to
interesting new examples of quasi-Hamiltonian $K$-manifolds, such as
the ``spinning $2n$-sphere'' for $K=\SU(n)$.  Secondly we construct a
universal (``master'') moduli space of parabolic bundles with
structure group $K$ over a marked Riemann surface.  The master moduli
space carries a natural action of a maximal torus of $K$ and a
torus-invariant stratification into manifolds, each of which has a
symplectic structure.  An essential ingredient in the construction is
the universal implosion.  Paradoxically, although the universal
implosion has no complex structure (it is the four-sphere for
$K=\SU(2)$), the master moduli space turns out to be a complex
algebraic variety.
\end{abstract}

\maketitle

\tableofcontents


\section{Introduction}\label{section;intro}

An important first step in understanding a symplectic manifold $M$
equipped with a Hamiltonian action of a compact group $K$ is a form of
abelianization.  This involves taking the inverse image under the
moment map $\Phi\colon M\to\k^*$ of a closed Weyl chamber
$\bar{\ca{C}}$ in the dual of the Lie algebra of a maximal torus $T$.
According to the Guillemin-Sternberg cross-section theorem the inverse
image $\Phi\inv(\ca{C})$ of the open chamber $\ca{C}$ is a smooth
symplectic submanifold of $M$ and carries a natural Hamiltonian
$T$-action, which encodes many of the properties of $M$.  However, the
preimage of the boundary facets of $\ca{C}$ is in general neither
smooth nor symplectic.  One obtains from $\Phi\inv(\bar{\ca{C}})$ a
symplectic variety (a stratified space with symplectic strata) by
"imploding" certain orbits in the boundary pieces.  This is the
subject of \cite{guillemin-jeffrey-sjamaar}.  As shown in that paper,
the implosion of the cotangent bundle $T^*K$ plays the role of a
universal implosion, from which all other implosions are derived in a
simple way.

In parallel to the theory of Hamiltonian $K$-manifolds, Alekseev,
Malkin, and Meinrenken introduced in
\cite{alekseev-malkin-meinrenken;lie-group} a notion of
quasi-Hamiltonian $K$-manifolds, which come equipped with moment maps
with values in the group $K$.  Quasi-Hamiltonian $K$-manifolds and
their moment maps share many of the features of the Hamiltonian ones,
such as reductions and cross-sections.  In this paper, we define the
notion of symplectic implosion of quasi-Hamiltonian $K$-manifolds,
where $K$ is a simply connected compact Lie group.  This
simultaneously generalizes the results of
\cite{guillemin-jeffrey-sjamaar} (where imploded cross-sections were
defined) and \cite{hurtubise-jeffrey;representations-weighted} (where
imploded cross-sections of particular quasi-Hamiltonian $K$-manifolds
were defined in the special case $K=\SU(2)$).  In the course of doing
this we define a universal implosion, obtained as the imploded
cross-section of the ``double'' $DK=K\times K$.

As an unexpected byproduct of this work we find new examples of
quasi-Hamiltonian $K$-manifolds, in particular the ``spinning
$2n$-sphere'', which generalizes the ``spinning four-sphere''
described in \cite{hurtubise-jeffrey;representations-weighted} and
\cite{alekseev-meinrenken-woodward;duistermaat-heckman}.

The motivation for developing this theory came from one particularly
important example, namely the space of representations of the
fundamental group of a punctured Riemann surface (obtained by
puncturing a closed Riemann surface at points $p_1$, $p_2$,\dots,
$p_n$) with values in a simply connected compact $K$.  Fixing the
conjugacy classes $C=(C_1,C_2,\dots,C_n)$ of the image in $K$ of small
loops around the punctures, one has symplectic moduli spaces
$M(\Sigma,C)$, about which some quite beautiful results are known, in
particular the theorem of Mehta and Sehadri \cite{mehta-seshadri}
saying that $M(\Sigma,C)$ is diffeomorphic to a space of parabolic
vector bundles over the closed Riemann surface.  One would like to
have a similar result about the space of all representations.  There
is no hope, however, of it being symplectic, as in many cases it is
not even even-dimensional.

We construct a ``master moduli space'' $M$ which is in some ways the
next best thing.  Let $\alcove\subset\t$ be the fundamental alcove of
$K$.  The image of $\bar{\alcove}$ under the exponential map
intersects every conjugacy class in $K$ in exactly one point.  The
space $M$ is a stratified space with symplectic strata equipped with a
Hamiltonian $T^n$-action and moment map
$M\to\bigl(\exp\bar{\alcove}\bigr)^n$.  The reduction of $M$ at the
$n$-tuple $C$ is then the parabolic moduli space $M(\Sigma,C)$ with
conjugacy classes $C$ at the punctures.

It is reasonable to hope that there is a version of the Mehta-Seshadri
theorem for the master moduli space $M$, and this is indeed the case:
for $K=\SU(n)$ one has a complex moduli space $\ca{M}$ of ``framed
parabolic sheaves'', and a homeomorphism of $M$ with $\ca{M}$.  The
parabolic moduli spaces $M(\Sigma,C)$ are GIT quotients of $\ca{M}$
under the action of the complexification of the group $T^n$.  The
complex algebraic variety $\ca{M}$ parametrizes sheaves on the Riemann
surface incorporating as extra data a parabolic structure at the
punctures and a framing of the successive quotients of the natural
flag attached to the parabolic structure.  The complexification of
$T^n$ acts transitively on these framings.  A surprising feature is
that the construction includes all possible types of flags, ranging
from the generic, maximal flag to the minimal, trivial one.  The type
of the flag for $M(\Sigma,C)$ at $p_i$ is determined by the face of
$\alcove$ to which $C_i$ belongs.  This complex picture is not treated
in this paper, but will be the subject of
\cite{hurtubise-jeffrey-sjamaar;moduli-framed}.

The layout of this paper is as follows.  In Section
\ref{section;hamiltonian} we recall the definition and basic
properties of quasi-Hamiltonian $K$-manifolds.  In Section
\ref{section;implosion} we recall the definition of imploded
cross-sections in the Hamiltonian case and extend it to the
quasi-Hamiltonian case.  In Section \ref{section;double} we define the
``universal implosion'' (the implosion of the double $DK$), after
reviewing the corresponding object in the Hamiltonian category (the
implosion of $T^*K$).  In Section \ref{section;projective} we describe
the universal imploded cross-section as the orbit of a toric variety
embedded in a representation of $K$.  In Section
\ref{section;parabolic} we descibe the construction of the master
moduli space via implosion.  In Appendix \ref{section;adjoint} we
summarize some basic facts concerning the conjugation action of a
compact Lie group and in Appendix \ref{section;sphere} we construct
the spinning $2n$-sphere.  A notation index can be found at the end of
the paper.

\section{Hamiltonian and quasi-Hamiltonian manifolds}
\label{section;hamiltonian}

\subsection*{Group-valued moment maps}

Let $K$\glossary{K@$K$, compact simply connected Lie group} be a
compact connected Lie group with Lie algebra $\k=T_1K$.  Recall that a
\emph{Hamiltonian $K$-manifold} is a symplectic manifold
$(M,\omega)$\glossary{M@$M$, (quasi-)Hamiltonian
$K$-manifold}\glossary{omega@$\omega$, (quasi-)symplectic form}
equipped with an action of $K$ which satisfies
\begin{equation}\label{equation;moment}
\iota(\xi_M)\omega=d\langle\Phi,\xi\rangle
\end{equation}
for an $\Ad^*$-equivariant map $\Phi\colon
M\to\k^*$\glossary{Phi@$\Phi$, moment map}, called the \emph{moment
map} for the action.  Here $\xi_M$ is the vector field generated by
$\xi\in\k$.  Alekseev et al.\
\cite{alekseev-malkin-meinrenken;lie-group} have extended this notion
in the following manner.  Recall that the Maurer-Cartan forms
$\theta_L$, $\theta_R\in\Omega^1(K,\k)$\glossary{theta@$\theta_L$,
$\theta_R$, Maurer-Cartan forms} are defined by
$\theta_{L,g}(L(g)_*\xi)=\xi$ and $\theta_{R,g}(R(g)_*\xi)=\xi$ for
$\xi\in\k$.  Here $L(g)$\glossary{Lg@$L(g)$, left multiplication by
$g$} denotes left multiplication and $R(g)$\glossary{Rg@$R(g)$, right
multiplication by $g$} right multiplication by $g$.  Fix a
$K$-invariant inner product $({\cdot},{\cdot})$ on $\k$ and let
$\chi=\frac1{12}(\theta_L,[\theta_L,\theta_L])\glossary{chi@$\chi$,
bi-invariant three-form} =\frac1{12}(\theta_R,[\theta_R,\theta_R])$ be
the corresponding bi-invariant three-form on $K$.

\begin{definition}\label{definition;quasi}
A \emph{quasi-Hamiltonian} (or \emph{group-valued Hamiltonian})
\emph{$K$-manifold} is a smooth $K$-manifold $M$\glossary{M@$M$,
(quasi-)Hamiltonian $K$-manifold} equipped with a $K$-invariant
two-form $\omega$\glossary{omega@$\omega$, (quasi-)symplectic form}
and an $\Ad$-equi\-var\-i\-ant\glossary{Ad@$\Ad$, adjoint action on
$K$ or $\k$} map $\Phi\colon M\to K$\glossary{Phi@$\Phi$, moment map},
called the \emph{(group-valued) moment map}, satisfying the following
properties:
\begin{enumerate}
\item\label{item;chi}
$d\omega=-\Phi^*\chi$;
\item\label{item;ker}
$\ker\omega_x=\{\,\xi_M(x)\mid\xi\in\ker(\Ad\Phi(x)+1)\,\}$ for all
$x\in M$;
\item\label{item;theta}
$\iota(\xi_M)\omega=\frac12\Phi^*(\theta_L+\theta_R,\xi)$ for all
$\xi\in\k$.
\end{enumerate}
\end{definition}

If $K$ is a torus, the two-form is symplectic and $\Phi$ is a moment
map in the sense of \cite{mcduff;moment-map}.  For nonabelian $K$ the
two-form is usually neither nondegenerate nor closed.  (But see
Appendix \ref{section;sphere} for a quasi-Hamiltonian
$\SU(n)$-structure on $\C\P^{n-1}$ whose two-form is symplectic.)
Axiom \eqref{item;chi} replaces the closedness of $\omega$, axiom
\eqref{item;ker}, often referred to as the \emph{minimal degeneracy}
axiom, replaces the nondegeneracy of $\omega$ and axiom
\eqref{item;theta} replaces the moment map condition
\eqref{equation;moment}.

Group-valued moment maps have rather limited functorial properties.  A
quasi-Hamiltonian action of $K$ usually does not restrict to
quasi-Hamiltonian actions of its subgroups.  For instance, if a
product $K_1\times K_2$ acts in a quasi-Hamiltonian fashion on
$(M,\omega)$, neither of the factors acts in a quasi-Hamiltonian
fashion (unless the other factor happens to be abelian).  However, the
following properties are straightforward consequences of the
definition.

\begin{lemma}\label{lemma;quasi}
Let $(M,\omega,\Phi)$ be a quasi-Hamiltonian $K$-manifold.
\begin{enumerate}
\item\label{item;centre}
For any $g$ in the centre of $K$, $L(g)\circ\Phi=R(g)\circ\Phi$ is a
moment map for the $K$-action on $(M,\omega)$.
\item\label{item;opposite}
$(M,-\omega,i\circ\Phi)$ is a quasi-Hamiltonian $K$-manifold, where
$i\colon K\to K$ denotes inversion.\glossary{i, inversion $g\mapsto
g\inv$}
\item\label{item;product}
If $K=K_1\times K_2$ where $K_2$ is a torus, then
$\pi_1\circ\Phi\colon M\to K_1$ is a moment map for the $K_1$-action
on $(M,\omega)$.
\item\label{item;cover}
Let $\pi_K\colon\tilde{K}\to K$ be a covering homomorphism and let
$\pi_M\colon\tilde{M}\to M$ be the induced covering, as in the
pullback diagram
$$
\xymatrix{\tilde{M}\ar@{->}[r]^-{\tilde{\Phi}}\ar@{->}[d]_{\pi_M}
&\tilde{K}\ar@{->}[d]^{\pi_K}\\M\ar@{->}[r]_-{\Phi}&**[r]K.}
$$
Then there is a unique $\tilde{K}$-action on $\tilde{M}$ such that
$\tilde{\Phi}$ is equivariant and $\pi_M$ is $\pi_K$-equivariant.  Let
$\tilde{\omega}=\pi_M^*\omega$.  Then $\tilde{\Phi}$ is a moment map
for the $\tilde{K}$-action on $(\tilde{M},\tilde{\omega})$.
\item\label{item;central}
Let $K_1$ be a closed central subgroup of $K$ which acts trivially on
$M$.  Let $K_2=K/K_1$ and let $\Phi_2\colon M\to K_2$ be the
composition of $\Phi$ with the quotient map $K\to K_2$.  Then the
triple $(M,\omega,\Phi_2)$ is a quasi-Hamiltonian $K_2$-manifold.
\end{enumerate}
\end{lemma}

The triple $(M,-\omega,i\circ\Phi)$ in \eqref{item;opposite} is the
quasi-Hamiltonian $K$-manifold \emph{opposite} to $M$.  We shall
frequently denote it by $M^-$.

If the moment map $\Phi$ in \eqref{item;cover} lifts to a map
$\Phi'\colon M\to\tilde{K}$, then $\tilde{M}$ is simply a disjoint
union of copies of $M$ and $\Phi'$ is the restriction of
$\tilde{\Phi}$ to one copy.



\subsection*{Fusion}

A sort of functoriality for group-valued moment maps holds for
restriction to a \emph{diagonal} subgroup.

\begin{theorem}[{\cite[Theorem
6.1]{alekseev-malkin-meinrenken;lie-group}}]\label{theorem;fusion}
Let $(M,\omega,\Phi)$ be a quasi-Hamiltonian $K\times K\times
H$-manifold, with moment map
$\Phi=\Phi_1\times\Phi_2\times\Phi_3\colon M\to K\times K\times H$.
Let $K\times H$ act on $M$ via the embedding $(k,h)\mapsto(k,k,h)$.
Then $M$ furnished with the two-form
$\omega+\frac12(\Phi_1^*\theta_L,\Phi_2^*\theta_R)$ and the moment map
$\Phi_1\Phi_2\times\Phi_3\colon M\to K\times H$ is a quasi-Hamiltonian
$K\times H$-manifold.
\end{theorem}

This restriction process from $K\times K\times H$ to $K\times H$ is
called \emph{internal fusion}.  The most important class of examples
is that of a Cartesian product $M_1\times M_2$ of a quasi-Hamiltonian
$K$-manifold $(M_1,\omega_1,\Phi_1)$ and a quasi-Hamiltonian $K\times
H$-manifold $(M_2,\omega_2,\Phi_2\times\Phi_3)$.  This is a
quasi-Hamiltonian $K\times K\times H$-manifold in an obvious way, and
fusing the two copies of $K$ gives rise to a quasi-Hamiltonian
$K\times H$-manifold called the \emph{fusion product} $M_1\fusion
M_2$.\glossary{..@$\fusion$, fusion product}

\subsection*{Exponentiation and linearization}

Let $(M,\omega_0,\Phi_0)$ be a Hamiltonian $K$-manifold in the
ordinary sense.  Let us use the isomorphism $\k\to\k^*$ given by the
inner product to identify $\k^*$ with $\k$, so that we can regard
$\Phi_0$ as a map into $\k$.  The process of \emph{exponentiation}
alters the two-form and the moment map (but not the action) on $M$,
namely into $\omega=\omega_0+\Phi_0^*\varpi$ and
$\Phi=\exp\circ\Phi_0$.  Here $\varpi\in\Omega^2(\k)$ is the $\Ad
K$-invariant primitive of $-\exp^*\chi\in\Omega^3(\k)$ given by
setting $\exp_s\lambda=\exp(s\lambda)$ and
$$
\varpi=\frac12\int_0^1\biggl(\exp_s^*\theta_R,\frac{\partial}{\partial
s}\exp_s^*\theta_R\biggr)\,ds.
$$
A calculation using
\begin{equation}\label{equation;maurer}
(\exp^*\theta_L)_\lambda=\frac{1-e^{-\ad\lambda}}{\ad\lambda}
\qquad\text{and}\qquad
(\exp^*\theta_R)_\lambda=\frac{e^{\ad\lambda}-1}{\ad\lambda}
\end{equation}
yields
\begin{equation}\label{equation;varpi}
\varpi_\lambda(\xi_1,\xi_2)=\int_0^1\biggl(\frac{1-\cosh
s\ad\lambda}{\ad\lambda}\xi_1,\xi_2\biggr)\,ds
=\biggl(\frac{\ad\lambda-\sinh\ad\lambda}{(\ad\lambda)^2}\xi_1,
\xi_2\biggr)
\end{equation}
for $\lambda$, $\xi_1$, $\xi_2\in\k$.  The triple $(M,\omega,\Phi)$
satisfies axioms \eqref{item;chi} and \eqref{item;theta} of Definition
\ref{definition;quasi}.  If in addition all points in $\Phi_0(M)$ are
regular for the exponential map, axiom \eqref{item;ker} is also
satisfied and so $(M,\omega,\Phi)$ is a quasi-Hamiltonian
$K$-manifold.  The reverse of exponentiation is \emph{linearization}.
Let $(M,\omega,\Phi)$ be a quasi-Ham\-ilton\-ian $K$-manifold.
Suppose there exists an $\Ad$-invariant open $U$ in $\k$ such that
$\exp\colon U\to K$ is a diffeomorphism onto an open subset containing
$\Phi(M)$ (with inverse denoted by $\log\colon\exp U\to U$).  The
linearization of $M$ is the Hamiltonian $K$-manifold
$(M,\omega_0,\Phi_0)$, where $\Phi_0=\log\circ\Phi$ and
$\omega_0=\omega-\Phi_0^*\varpi$.  (See
\cite[\S~3.3]{alekseev-malkin-meinrenken;lie-group}.)

\subsection*{Quasi-symplectic quotients}

It is shown in \cite[\S~8]{alekseev-malkin-meinrenken;lie-group} that
the category of quasi-Hamiltonian $K$-manifolds is equivalent to a
subcategory of the category of infinite-dimensional symplectic
manifolds with Hamiltonian actions (in the usual sense) of the
\emph{loop group} $LK$.  It should therefore come as no surprise that
many constructions in the category of Hamiltonian manifolds have
parallels in the quasi-Hamiltonian world.  In particular there are
analogues of symplectic reduction and of symplectic implosion.  We
review quasi-symplectic reduction below; quasi-symplectic implosion is
the topic of Section \ref{section;implosion}.

Let $(M,\omega,\Phi)$ be a quasi-Hamiltonian $K$-manifold.  Assume
that $K=K_1\times K_2$ where $K_2$ is a torus and let $\Phi_1\colon
M\to K_1$ and $\Phi_2\colon M\to K_2$ be the components of $\Phi$.
Let $g\in K$.  Because $\Phi$ is equivariant, the centralizer
$(K_1)_g$ acts on the fibre $\Phi_1\inv(g)$.  The
\emph{quasi-symplectic quotient} or \emph{quasi-Hamiltonian reduced
space} at $g$ is the topological space
$$
M\quot[g]K_1=\Phi_1\inv(g)/(K_1)_g.\glossary{M@$M\quot[g]K$,
(quasi-)Hamiltonian quotient}
$$
The subscript is usually omitted when $g=1$.  In good cases this
quotient is a symplectic orbifold.

\begin{theorem}%
[{\cite[Theorem 5.1]{alekseev-malkin-meinrenken;lie-group}}]
\label{theorem;quotient}
Suppose that $g$ is a regular value of $\Phi_1$.  Then the centralizer
$(K_1)_g$ acts locally freely on the submanifold $\Phi_1\inv(g)$.  The
restriction of $\omega$ to $\Phi_1\inv(g)$ is closed and
$(K_1)_g$-basic.  The form $\omega_g$ on the orbifold $M\quot[g]K_1$
induced by $\omega$ is nondegenerate.  The map $M\quot[g]K_1\to K_2$
induced by $\Phi_2$ is a moment map for the induced $K_2$-action on
$M\quot[g]K_1$.
\end{theorem}

(\cite{alekseev-malkin-meinrenken;lie-group} also covers the case
where the second factor $K_2$ is nonabelian.  Then $M\quot[g]K$ is not
symplectic, but a quasi-Hamiltonian $K_2$-orbifold.)

In the singular case the quotient stratifies into symplectic manifolds
according to orbit type.  Let $H$ be a subgroup of $(K_1)_g$.  Recall
that the stratum of orbit type $H$ in $M$ (with respect to the
$(K_1)_g$-action) is the $(K_1)_g$-invariant submanifold $M_{(H)}$
consisting of all points $x$ such that the stabilizer
$(K_1)_g\cap(K_1)_x$ is conjugate to $H$.  Put $Z=\Phi_1\inv(g)$ and
$Z_{(H)}=Z\cap M_{(H)}$.  Let $\{\,Z_i\mid i\in I\,\}$ be the
collection of connected components of all subspaces $Z_{(H)}$, where
$(H)$ ranges over all conjugacy classes of subgroups of $(K_1)_g$.
Set-theoretically the quotient is a disjoint union,
\begin{equation}\label{equation;orbit-type}
M\quot[g]K_1=\coprod_{i\in I}Z_i/(K_1)_g.
\end{equation}

\begin{theorem}\label{theorem;singular-quotient}
Let $g\in K_1$ be arbitrary.  Each of the subsets $Z_i$ is a
submanifold of $M$ and the restriction of $\omega$ to $Z_i$ is closed
and $(K_1)_g$-basic.  The orbit space $Z_i/(K_1)_g$ is a manifold and
the form induced by $\omega$ on $Z_i/(K_1)_g$ is nondegenerate.  The
decomposition \eqref{equation;orbit-type} is a locally normally
trivial stratification of $M\quot[g]K_1$.  The stratification is
$K_2$-invariant and the continuous map $\bar{\Phi}_2\colon
M\quot[g]K_1\to K_2$ induced by $\Phi_2$ restricts to a moment map for
the $K_2$-action on each stratum.
\end{theorem}

\begin{proof}
By the quasi-Hamiltonian shifting trick (see
\cite[Remark~6.2]{alekseev-malkin-meinrenken;lie-group}) we may assume
that $g=1$.  Choose an invariant open neighbourhood $U$ of $0\in\k$
such that $\exp$ is a diffeomorphism from $U$ onto its image and let
$\log\colon\exp U\to U$ be the inverse.  Replace $M$ with the (not
necessarily connected) open subset $\Phi\inv(\exp U)$.  Clearly this
does not affect the quotient $M\quot K$.  Let $(M,\omega_0,\Phi_0)$ be
the linearization of $(M,\omega,\Phi)$.  Observe that
$\Phi\inv(1)=\Phi_0\inv(0)$.  Hence $\omega=\omega_0$ on every
submanifold of $M$ contained in $\Phi\inv(1)$.  The theorem follows
therefore from the Hamiltonian case, where we appeal to \cite[Theorem
2.1]{sjamaar-lerman;stratified}.
\end{proof}

(Presumably this result too generalizes to nonabelian $K_2$.  However,
because the linearization $(M,\omega_0,\Phi_0)$ is not a Hamiltonian
$K$-manifold in that case, our proof is valid only for abelian $K_2$.)
We call the space $M\quot[g]K_1$ a \emph{stratified Hamiltonian
$K_2$-space} and refer to the map $\bar{\Phi}_2$ as the moment map for
the $K_2$-action.

\section{Imploded cross-sections}\label{section;implosion}

\subsection*{The Hamiltonian case}

The first goal of this paper is to develop a quasi-Ham\-ilton\-ian
analogue of symplectic implosion.  We start by reviewing this notion
in the Hamiltonian case, referring to
\cite[\S~2]{guillemin-jeffrey-sjamaar} for details.  Symplectic
implosion is an ``abelianization functor'', a crude operation that
transmutes a Hamiltonian $K$-manifold into a Hamiltonian $T$-space
(where $T$ is a maximal torus of $K$) retaining some of the relevant
features of the original manifold, but at the expense of producing
singularities.  However, the singular set breaks up into smooth
symplectic manifolds in a nice way.

Fix a maximal torus $T$\glossary{T@$T$, maximal torus of $K$} of $K$
and an open chamber $\ca{C}$\glossary{C@$\ca{C}$, open chamber in
$\t^*$} in $\t^*$, the dual of $\t=\Lie T$.  The closed chamber
$\bar{\ca{C}}$\glossary{C@$\bar{\ca{C}}$, closure of $\ca{C}$} is a
polyhedral cone, which is the disjoint union of $2^r$ relatively open
faces, where $r$ is the rank of the commutator subgroup $[K,K]$.  We
define a partial order on the faces by putting
$\sigma\le\tau$\glossary{sigma@$\sigma$, face of chamber or
alcove}\glossary{.@$\le$, inclusion of faces} if
$\sigma\subset\bar{\tau}$.

Let $(M,\omega,\Phi)$ be a connected Hamiltonian $K$-manifold.  The
\emph{principal face} $\sigma\prin$\glossary{sigma@$\sigma\prin$,
principal face} is the smallest face $\sigma$ of $\ca{C}$ such that
the Kirwan polytope $\Phi(M)\cap\bar{\ca{C}}$ is contained in the
closure of $\sigma$.  In many cases $\sigma\prin=\ca{C}$.  The
\emph{cross-section} of $M$ is $\Phi^{-1}(\sigma\prin)$.  This is a
$T$-invariant connected symplectic submanifold of $M$.  The torus
action on the cross-section is Hamiltonian with moment map equal to
the restriction of the $K$-moment map, and the $K$-invariant subset
$K\Phi^{-1}(\sigma\prin)$ is open and dense in $M$.

The \emph{imploded cross-section} is a ``completion'' of the
cross-section to a stratified space with symplectic strata.  It is
obtained by taking the preimage of the closed chamber,
$\crosschamber$, which stratifies into smooth manifolds in a natural
way, and by quotienting out the null-foliation of the form $\omega$ on
each stratum.  To wit, declare two points $m_1$ and $m_2$ in
$\crosschamber$ to be \emph{equivalent} if there exists $k$ in the
commutator group $\bigl[K_{\Phi(m_1)},K_{\Phi(m_1)}\bigr]$ such that
$m_2=km_1$.  (Here $K_\xi$ denotes the centralizer of $\xi\in\k^*$.)
Then the imploded cross-section is the quotient space
$$
M\impl=\crosschamber/{\sim},\glossary{M@$M\impl$, imploded
cross-section}\glossary{impl, implosion}\glossary{..@$\sim$,
equivalence relation defining $M\impl$}
$$
equipped with the quotient topology.  Set-theoretically it is a
disjoint union
\begin{equation}\label{equation;decomposition}
M\impl=\coprod_{\sigma\le\ca{C}}
\Phi\inv(\sigma)\big/[K_\sigma,K_\sigma]
\end{equation}
over the faces $\sigma$ of $\ca{C}$.  Here
$K_\sigma$\glossary{K@$K_\sigma$, centralizer of a face} is the
centralizer of (any point in) the face $\sigma$.  The pieces in the
decomposition \eqref{equation;decomposition} are usually not
manifolds.  However, the decomposition can be refined into a
stratification of $M\impl$ with symplectic strata.  There is a unique
open stratum, which is dense in $M\impl$ and symplectomorphic to the
cross-section $\Phi^{-1}(\sigma\prin)$.  (See \cite[Theorems 2.10 and
5.10]{guillemin-jeffrey-sjamaar}.)  These facts are most easily
established by using slices for the coadjoint action on $\k^*$.  The
\emph{natural slice} at a face $\sigma$ is the subset of $\k_\sigma^*$
given by
$$
\eu{S}_\sigma=\Ad^*(K_\sigma)\star\sigma,\glossary{S@$\eu{S}_\sigma$,
slice at face $\sigma$}
$$
where $\star\sigma$\glossary{star@$\star\sigma$, star of face
$\sigma$} denotes the open star $\bigcup_{\tau\ge\sigma}\tau$ of
$\sigma$.  ($\eu{S}_\sigma$ is a slice for the coadjoint action in the
sense that $K\eu S_\sigma$ is open in $\k^*_\sigma$ and
$K$-equivariantly diffeomorphic to the associated bundle
$K\times^{K_\sigma}\eu S_\sigma$.  Note that $\eu{S}_\sigma=\sigma$ if
$\sigma=\ca{C}$.)  The set $M_\sigma=\Phi\inv(\eu
S_\sigma)$\glossary{M@$M_\sigma$, cross-section over $\sigma$} is a
$K_\sigma$-invariant connected symplectic submanifold of $M$ and
$\Phi\inv(\sigma)/[K_\sigma,K_\sigma]$ is the symplectic quotient of
$M_\sigma$ with respect to the action of $[K_\sigma,K_\sigma]$.
Thence, by appealing to \cite[Theorem 2.1]{sjamaar-lerman;stratified}
we obtain a refinement of the decomposition
\eqref{equation;decomposition} into symplectic strata.

The \emph{imploded moment map}\glossary{Phi@$\Phi\impl$, imploded
moment map} is the continuous map $\Phi\impl\colon M\impl\to\t^*$
induced by $\Phi$.  The imploded cross-section is an abelianization of
$M$ in the sense that it carries a natural $T$-action which preserves
the stratification.  The restriction of $\Phi\impl$ to each stratum is
a moment map for the $T$-action; the image of $\Phi\impl$ is equal to
$\Phi(M)\cap\bar{\ca{C}}$, the Kirwan polytope of $M$; and the
symplectic quotients of $M\impl$ with respect to the $T$-action are
identical to the symplectic quotients of $M$ with respect to the
$K$-action.  (See \cite[Theorem 3.4]{guillemin-jeffrey-sjamaar}.)

\subsection*{The quasi-Hamiltonian case}

Now let $(M,\omega,\Phi)$ be a connected quasi-Hamilton\-ian
$K$-manifold.  Proceeding by analogy with the Hamiltonian case we
shall define the imploded cross-section of $M$ and prove that it
stratifies into quasi-Hamiltonian $T$-manifolds.  To obtain good
results we shall assume for the remainder of the paper that
$K$\glossary{K@$K$, compact simply connected Lie group} is
\emph{simply connected}.

Let $\ca{C}\spcheck$\glossary{C@$\ca{C}\spcheck$, dual open chamber in
$\t$} be the chamber in $\t$ dual to $\ca{C}$ and let
$\alcove$\glossary{A@$\alcove$, open alcove in $\t$} be the unique
(open) alcove contained in $\ca{C}\spcheck$ such that
$0\in\bar{\alcove}$.\glossary{A@$\bar{\alcove}$, closure of $\alcove$}
The exponential map induces a homeomorphism from
$\bar{\alcove}/\pi_1(K)=\bar{\alcove}$ to $T/W\cong K/\Ad K$, the
space of conjugacy classes in $K$.  (Cf.\ \cite[Ch.~9, \S~5.2,
Corollaire 1]{bourbaki;groupes-algebres}.)  Here
$W=N(T)/T$\glossary{W@$W$, Weyl group} denotes the Weyl group of
$(K,T)$.  In particular the restriction of the exponential map to
$\bar{\alcove}$ is injective.  For $g$ in $K$ let $K_g=\{\,h\in K\mid
hg=gh\,\}$\glossary{K@$K_g$, centralizer of $g\in K$} denote the
centralizer of $g$.  For points $m_1$, $m_2\in\cross$ define $m_1\sim
m_2$ if $m_2=km_1$ for some
$k\in\bigl[K_{\Phi(m_1)},K_{\Phi(m_1)}\bigr]$.  Since $\Phi$ is
$\Ad$-equivariant, $m_1\sim m_2$ implies $\Phi(m_1)=\Phi(m_2)$, and
hence $\sim$ is an equivalence relation.

\begin{definition}\label{definition;implosion}
The \emph{imploded cross-section} of $M$ is the quotient space
$M\impl=\cross/{\sim}$,\glossary{M@$M\impl$, imploded
cross-section}\glossary{impl, implosion}\glossary{..@$\sim$,
equivalence relation defining $M\impl$} equipped with the quotient
topology.  The \emph{imploded moment map}
$\Phi\impl$\glossary{Phi@$\Phi\impl$, imploded moment map} is the
continuous map $M\impl\to T$ induced by $\Phi$.
\end{definition}

Imitating \cite[Lemma 2.3]{guillemin-jeffrey-sjamaar} one proves
easily that the quotient map from $\cross$ to $M\impl$ is proper and
that $M\impl$ is Hausdorff, locally compact and second countable.  The
action of $T$ preserves $\cross$ and descends to a continuous action
on $M\impl$.  The image of the imploded moment map is
$\Phi(M)\cap\exp\bar{\alcove}$, the moment polytope of $M$.

The imploded cross-section is seldom a manifold, but we shall see that
it partitions into manifolds, each of which carries a symplectic
structure.  The following discussion draws on some standard facts
about the conjugation action reviewed in Appendix
\ref{section;adjoint}.  Lemma \ref{lemma;borel} shows that we have an
analogue of the decomposition \eqref{equation;decomposition} into
orbit spaces,
\begin{equation}\label{equation;quasi-decomposition}
M\impl=\coprod_{\sigma\le\alcove}
\Phi\inv(\exp\sigma)\big/[K_\sigma,K_\sigma],\glossary{.@$\le$,
inclusion of faces}
\end{equation}
where this time $\sigma$ ranges over the faces of the alcove $\alcove$
and $K_\sigma$ stands for the centralizer of $\exp\sigma$.  Let us
write $X_\sigma
=\Phi\inv(\exp\sigma)\big/[K_\sigma,K_\sigma]$.\glossary{X@$X_\sigma$,
piece $\Phi\inv(\exp\sigma)\big/[K_\sigma,K_\sigma]$ of $M\impl$} Each
of the $X_\sigma$ is invariant under the $T$-action and the imploded
moment map carries $X_\sigma$ into $\exp\sigma$.

Fix a face $\sigma\le\alcove$ and consider its open star
$\star\sigma=\bigcup_{\tau\ge\sigma}\tau$\glossary{star@$\star\sigma$,
star of face $\sigma$}.  By Proposition \ref{proposition;slice} the
set $\eu{S}_\sigma=\Ad(K_\sigma)\exp(\star\sigma)\subset
K_\sigma$\glossary{S@$\eu{S}_\sigma$, slice at face $\sigma$} is a
slice for the conjugation action of $K$ on itself.  Equivariant maps
are transverse to slices, so the $K_\sigma$-invariant subset
$M_\sigma=\Phi\inv(\eu{S}_\sigma)$\glossary{M@$M_\sigma$,
cross-section over $\sigma$} is a submanifold of $M$, called the
\emph{cross-section} of $M$ over $\sigma$.  As in the Hamiltonian
case, the \emph{principal face}
$\sigma\prin$\glossary{sigma@$\sigma\prin$, principal face} for $M$ is
defined as the smallest $\sigma\le\alcove$ such that the moment
polytope $\Phi(M)\cap\bar{\alcove}$ is contained in the closure of
$\sigma$.  The \emph{principal cross-section} is $M_{\sigma\prin}$.

\begin{theorem}[quasi-Hamiltonian cross-sections]
\label{theorem;cross}
Let $\sigma$ be a face of $\alcove$.
\begin{enumerate}
\item\label{item;cross}
$M_\sigma$ is a connected quasi-Hamiltonian $K_\sigma$-manifold with
moment map $\Phi_\sigma=\Phi|_{M_\sigma}$\glossary{Phi@$\Phi_\sigma$,
moment map on $M_\sigma$}.
\item\label{item;dense}
The action $K\times M_\sigma\to M$ induces a diffeomorphism
$K\times^{K_\sigma}M_\sigma\to KM_\sigma$.  If $M_\sigma$ is nonempty,
then $KM_\sigma$ is dense in $M$.
\item\label{item;principal}
The commutator subgroup of $K_{\sigma\prin}$ acts trivially on the
principal cross-section $M_{\sigma\prin}$.  Hence $M_{\sigma\prin}$ is
symplectic.
\end{enumerate}
\end{theorem}

\begin{proof}
See \cite[\S~7]{alekseev-malkin-meinrenken;lie-group} for a proof that
$M_\sigma$ is a quasi-Hamiltonian $K_\sigma$-manifold.

The diffeomorphism $K\times^{K_\sigma}M_\sigma\cong KM_\sigma$ follows
from the slice theorem, Proposition \ref{proposition;slice}.  The
connectedness of $M_\sigma$ and denseness of $KM_\sigma$ are proved as
in the Hamiltonian case, \cite[\S~3]{lerman;nonabelian-convexity},
using the simple connectivity of the homogeneous space $K/K_\sigma$.

Put $\sigma=\sigma\prin$.  Then
$\Phi(M)\cap\exp\bar{\alcove}\subset\bar{\sigma}$, so
$M_\sigma=\Phi\inv(\exp\sigma)$.  In particular,
$K_{\Phi(x)}=K_\sigma$ for all $x\in M_\sigma$.  Now fix $x\in
M_\sigma$ and $v\in T_xM_\sigma$.  Select a curve $x_t$ in $M_\sigma$
with $x_0=x$ and $x'_0=v$.  Let $\xi\in\k_\sigma$.  Since
$\k_\sigma=\ker(\Ad(\Phi(x_t))-1)$ we have $\Ad(\Phi(x_t))\xi=\xi$ for
all $t$.  Differentiating this identity and setting $t=0$ yields
$$
[\xi,(\Phi_\sigma^*\theta_{L,\sigma})_x(v)]
=[\xi,(\Phi^*\theta_L)_x(v)]
=\bigl[\xi,L\bigl(\Phi(x)\inv\bigr)_*T_x\Phi(v)\bigr]=0.
$$
Here $\theta_{L,\sigma}$ denotes the Maurer-Cartan form on $K_\sigma$.
This means that the image of $(\Phi_\sigma^*\theta_{L,\sigma})_x\colon
T_xM_\sigma\to\k_\sigma$ is contained in the centre of $\k_\sigma$.
According to \cite[Proposition
4.1.3]{alekseev-malkin-meinrenken;lie-group} this image is equal to
the orthogonal complement $\k_x^\perp$ of $\k_x$ in $\k_\sigma$.  In
other words, $\k_x^\perp\subset\z(\k_\sigma)$, which implies
$[\k_\sigma,\k_\sigma]\subset\k_x$.  Therefore $[K_\sigma,K_\sigma]$
acts trivially on $M_\sigma$.  This shows that $M_\sigma$ is a
quasi-Hamiltonian manifold for the torus
$K_\sigma/[K_\sigma,K_\sigma]$ and in particular is symplectic.
\end{proof}

\begin{remark}\label{remark;product}
Suppose $K$ is a product $K_1\times K_2$ and $M$ is quasi-Hamiltonian
with moment map $\Phi=\Phi_1\times\Phi_2$.  Let
$\alcove=\alcove_1\times\alcove_2$ be the alcove of $K$ and let
$\sigma$ be a face of $\alcove_2$.  Let $\eu{S}_\sigma$ be the
corresponding slice in $K_2$, and let
$M_\sigma=\Phi_2\inv(\eu{S}_\sigma)$.  Then $M_\sigma$ is a connected
quasi-Hamiltonian $K_1\times(K_2)_\sigma$-manifold with moment map
$\Phi_\sigma=\Phi|_{M_\sigma}$.  This is proved by noting that
$K_1\times\eu{S}_\sigma$ is a slice for the adjoint action of $K$,
that $M_\sigma=\Phi\inv(K_1\times\eu{S}_\sigma)$, and by subsequently
applying Theorem \ref{theorem;cross}\eqref{item;cross}.
\end{remark}

As in \cite[Corollary 2.6]{guillemin-jeffrey-sjamaar} it follows from
Theorem \ref{theorem;cross} that the imploded cross-section contains a
dense copy of the principal cross-section.  We can thus view $M\impl$
as a ``completion'' of the symplectic manifold $M_{\sigma\prin}$.

\begin{corollary}
The restriction of the quotient map $\cross\to M\impl$ to the
principal cross-section is a homeomorphism onto its image.  The image
is open and dense in $M\impl$.  Hence $M\impl$ is connected.
\end{corollary}

We are now going to show that each of the pieces
$X_\sigma=\Phi\inv(\exp\sigma)\big/[K_\sigma,K_\sigma]$ of the
imploded cross-section can be written as a quasi-symplectic quotient
(up to a finite covering map).  Let $Z(K_\sigma)$ be the centre of
$K_\sigma$ and $Z(K_\sigma)^0$ its unit component.  Let $H_\sigma$ be
the affine subspace of $\t$ spanned by $\sigma$ and choose
$g_\sigma\in[K_\sigma,K_\sigma]\cap Z(K_\sigma)$ such that $\exp
H_\sigma\subset g_\sigma Z(K_\sigma)^0$, as in Lemma
\ref{lemma;face}\eqref{item;face-shift}.  By Lemma
\ref{lemma;quasi}\eqref{item;centre}
$\Phi_\sigma'=L(g_\sigma\inv)\circ\Phi_\sigma$ is a moment map for the
$K_\sigma$-action on $M_\sigma$.  The reason for shifting the moment
map is that it enables us to write $\Phi\inv(\exp\sigma)$ as the
inverse image of the torus $Z(K_\sigma)^0$.

\begin{lemma}\label{lemma;shift-slice}
\begin{enumerate}
\item\label{item;slice-centre}
$\exp\sigma=\eu{S}_\sigma\cap g_\sigma Z(K_\sigma)^0=\eu{S}_\sigma\cap
Z(K_\sigma)$.
\item\label{item;inverse-centre}
$\Phi\inv(\exp\sigma)=(\Phi_\sigma')\inv\bigl(Z(K_\sigma)^0\bigr)$.
\end{enumerate}
\end{lemma} 

\begin{proof}
The inclusion $\exp\sigma\subset\eu{S}_\sigma\cap g_\sigma
Z(K_\sigma)^0$ follows from the choice of $g_\sigma$ and the inclusion
$\eu{S}_\sigma\cap g_\sigma Z(K_\sigma)^0\subset\eu{S}_\sigma\cap
Z(K_\sigma)$ is obvious.  Now let $g\in\eu{S}_\sigma\cap Z(K_\sigma)$
and write $g=\Ad(k)\exp\xi$ with $k\in K_\sigma$ and
$\xi\in\star\sigma$.  The assumption $\Ad(k)\exp\xi\in Z(K_\sigma)$
implies $\Ad(k)\exp\xi\in T$, so there exists $w$ in the Weyl group
$N_{K_\sigma}(T)/T$ such that $\Ad(k)\exp\xi=w\exp\xi$.  Using
\eqref{equation;centre} and the fact that $w$ preserves $R_\sigma$,
the root system of $K_\sigma$, we see that $w\exp\xi\in Z(K_\sigma)$
implies $\alpha(\xi)\in2\pi i\Z$ for all $\alpha\in R_{\sigma,+}$.
Since $\xi\in\star\sigma$ we obtain $\alpha(\xi)=\alpha(\sigma)$ from
\eqref{equation;inequality-star}, and therefore $\xi\in\sigma$ by
\eqref{equation;inequality-face}.  Hence $\Ad(k)\exp\xi\in\exp\sigma$,
i.e.\ $\eu{S}_\sigma\cap Z(K_\sigma)\subset\exp\sigma$.

Let $x\in M_\sigma$.  Using \eqref{item;slice-centre} we see that
$\Phi_\sigma'(x)\in Z(K_\sigma)^0$ if and only if
$\Phi_\sigma(x)=g_\sigma\Phi_\sigma'(x)\in\exp\sigma$.
\end{proof}

Let us write $K_1=[K_\sigma,K_\sigma]$ and $K_2=Z(K_\sigma)^0$.
Recalling that $K_\sigma=K_1K_2$ we define the finite covering
\begin{equation}\label{equation;covering}
\tilde{K}_\sigma=K_1\times K_2.
\end{equation}
The covering map $\tilde{K}_\sigma\to K_\sigma$ is given by
multiplication and the covering group is the group
$\Gamma_\sigma=K_1\cap K_2$\glossary{Gamma@$\Gamma_\sigma$, covering
group $[K_\sigma,K_\sigma]\cap Z(K_\sigma)^0$} of Lemma
\ref{lemma;face}\eqref{item;face-shift}, embedded into
$\tilde{K}_\sigma$ by the antidiagonal map $g\mapsto(g,g\inv)$.
Pulling back this covering via the shifted moment map
$\Phi_\sigma'\colon M_\sigma\to K_\sigma$ yields a normal covering map
$\pi_\sigma\colon\tilde{M}_\sigma\to M_\sigma$ with the same covering
group $\Gamma_\sigma$.  By Lemma \ref{lemma;quasi}\eqref{item;cover}
$\tilde{M}_\sigma$ is a quasi-Hamiltonian $\tilde{K}_\sigma$-manifold
with a moment map
$\tilde{\Phi}_\sigma\colon\tilde{M}_\sigma\to\tilde{K}_\sigma$
obtained by lifting $\Phi_\sigma'$.  By Lemma
\ref{lemma;quasi}\eqref{item;product} the map
$\pi_1\circ\tilde{\Phi}_\sigma$ is a moment map for the $K_1$-action
on $\tilde{M}_\sigma$, where $\pi_1\colon\tilde{K}_\sigma\to K_1$ is
the projection.  Consider the quasi-symplectic quotient
$$
\tilde{X}_\sigma=\tilde{M}_\sigma\quot K_1
=(\pi_1\circ\tilde{\Phi}_\sigma)\inv(1)\big/K_1
=\tilde{\Phi}_\sigma\inv(K_2)\big/K_1,
$$
a space which, by Theorem \ref{theorem;singular-quotient}, stratifies
naturally into quasi-Hamiltonian $K_2$-manifolds.  The moment map is
the continuous map
\begin{equation}\label{equation;centre-moment}
\Psi_2\colon\tilde{X}_\sigma\to K_2
\end{equation}
induced by $\pi_2\circ\tilde{\Phi}_\sigma\colon\tilde{M}_\sigma\to
K_2$.  Recall that $X_\sigma$ denotes the piece
$\Phi\inv(\exp\sigma)\big/K_2$ in the decomposition
\eqref{equation;quasi-decomposition}.

\begin{lemma}\label{lemma;piece}
\begin{enumerate}
\item\label{item;preimage}
$\tilde{\Phi}_\sigma\inv(K_2)=\pi_\sigma\inv\Phi\inv(\exp\sigma)$.
\item\label{item;covering-action}
The $\Gamma_\sigma$-action on $\tilde{M}_\sigma$ descends to a free
action on $\tilde{X}_\sigma$.
\item\label{item;quotient-covering}
$\pi_\sigma$ induces a normal covering map
$\pi_\sigma\colon\tilde{X}_\sigma\to X_\sigma$ with covering group
$\Gamma_\sigma$.  The $K_2$-action on $\tilde{X}_\sigma$ commutes with
the $\Gamma_\sigma$-action.  The $K_2$-moment map
\eqref{equation;centre-moment} and the reduced symplectic forms on the
strata of $\tilde{X}_\sigma$ are $\Gamma_\sigma$-invariant.
\end{enumerate}
\end{lemma}

\begin{proof}
\eqref{item;preimage} follows from Lemma
\ref{lemma;shift-slice}\eqref{item;inverse-centre} and the fact that
$\tilde{\Phi}_\sigma$ is the lift of $\Phi_\sigma'$.

By definition $\tilde{M}_\sigma$ is the fibred product
\begin{equation}\label{equation;pullback}
\tilde{M}_\sigma=\{\,(x,k_1,k_2)\mid\text{$x\in M_\sigma$, $k_1\in
K_1$, $k_2\in K_2$, $\Phi_\sigma'(x)=k_1k_2$}\,\}.
\end{equation}
The action of $g\in\Gamma_\sigma$ is given by
$g\cdot(x,k_1,k_2)=(x,gk_1,g\inv k_2)$, whereas the action of $k\in
K_1$ is given by $k\cdot(x,k_1,k_2)=(kx,kk_1k\inv,k_2)$.  Therefore
$g\cdot(x,k_1,k_2)=k\cdot(x,k_1,k_2)$ implies $g=1$, so
$\Gamma_\sigma$ acts freely on the orbit space
$\tilde{X}_\sigma=\tilde{M}_\sigma/K_1$.

The first statement in \eqref{item;quotient-covering} follows
immediately from \eqref{item;preimage} and
\eqref{item;covering-action}.  Since $\Gamma_\sigma$ is central in
$\tilde{K}_\sigma$, its action on $\tilde{X}_\sigma$ commutes with
that of $K_2$ and $\Psi_2$ is $\Gamma_\sigma$-invariant.  The lift
$\tilde{\omega}_\sigma$ of the two-form $\omega_\sigma$ on $M_\sigma$
is $\tilde{K}_\sigma$-invariant, so the induced forms on the strata of
$\tilde{X}_\sigma$ are $\Gamma_\sigma$-invariant.
\end{proof}

Thus we have shown that each of the subspaces $X_\sigma\subset M\impl$
can be written, up to a finite covering map, as a quasi-symplectic
quotient.  Recall from Theorem \ref{theorem;singular-quotient} that
the covering space $\tilde{X}_\sigma$ stratifies according to
$K_1$-orbit type into symplectic manifolds.  Because $\Gamma_\sigma$
is central in $K_1$, this stratification is $\Gamma_\sigma$-invariant.
From Lemma \ref{lemma;piece}\eqref{item;quotient-covering} we now
obtain the following.

\begin{corollary}\label{corollary;stratification}
$X_\sigma$ stratifies into quasi-Hamiltonian $K_2$-manifolds.  The
strata are of the form $S/\Gamma_\sigma$, where $S$ runs over the
strata of $\tilde{X}_\sigma$.
\end{corollary}

The $K_2$-moment map sends $X_\sigma$ into
$g_\sigma\inv\exp\sigma\subset K_2$.  Moreover, the $K_2$-action on
$X_\sigma$ extends naturally to the maximal torus $T$ in the following
way.  First extend the $K_2$-action on $\tilde{X}_\sigma$ to a
quasi-Hamiltonian action of the maximal torus
$\tilde{T}_\sigma=(K_1\cap T)\times K_2$ of $\tilde{K}_\sigma$ by
letting the factor $K_1\cap T$ act trivially.  As moment map
$\tilde{\Psi}\colon\tilde{X}_\sigma\to\tilde{T}_\sigma$ we take
$\tilde{\Psi}=(\Psi_1,\Psi_2)$, where $\Psi_1\colon\tilde{X}_\sigma\to
K_1\cap T$ is the constant map sending $\tilde{X}_\sigma$ to
$g_\sigma\in K_1\cap Z(K_\sigma)\subset K_1\cap T$ and $\Psi_2$ is as
in \eqref{equation;centre-moment}.  Then $\tilde{\Psi}$ is
$\Gamma_\sigma$-equivariant and so induces a moment map $\Psi$ for the
action of $T=\tilde{T}_\sigma/\Gamma_\sigma$ on $X_\sigma$, as in the
diagram
$$
\xymatrix{\tilde{X}_\sigma\ar@{->}[r]^-{\tilde{\Psi}}
\ar@{->}[d]_{/\Gamma_\sigma}&\tilde{T}\ar@{->}[d]^{/\Gamma_\sigma}\\
X_\sigma\ar@{->}[r]_-{\Psi}&**[r]T.}
$$

\begin{lemma}\label{lemma;implosion}
$\Psi$ is equal to the restriction to $X_\sigma$ of the imploded
moment map $\Phi\impl$.
\end{lemma}

\begin{proof}
We need to show that the following diagram commutes:
\begin{equation}\label{equation;moment-commute}
\vcenter{\xymatrix{\tilde{\Phi}_\sigma\inv(K_2)\ar@{->}[r]^-{/K_1}
\ar@{->}[d]_{/\Gamma_\sigma}
&\tilde{X}_\sigma\ar@{->}[r]^-{\tilde{\Psi}}
&\tilde{T}_\sigma\ar@{->}[d]^{/\Gamma_\sigma}\\
\Phi\inv(\exp\sigma)\ar@{->}[r]_-{/K_1}
&X_\sigma\ar@{->}[r]_-{\Phi\impl}&**[r]T.}}
\end{equation}
The moment map on the fibred product \eqref{equation;pullback} is the
projection $\tilde{\Phi}_\sigma(x,k_1,k_2)=(k_1,k_2)$.  Therefore
$$
\tilde{\Phi}_\sigma\inv(K_2)=\{\,(x,1,k_2)\mid\text{$x\in M_\sigma$,
$k_2\in K_2$, $\Phi_\sigma'(x)=k_2$}\,\},
$$
and the composition of the two horizontal arrows at the top of
\eqref{equation;moment-commute} is the map
$(x,1,k_2)\mapsto(g_\sigma,\Phi_\sigma'(x))$.
Composing with the quotient map $\tilde{T}_\sigma\to T$, which is
given by multiplication, we get the map
\begin{equation}\label{equation;arrows}
(x,1,k_2)\mapsto g_\sigma\Phi_\sigma'(x)=\Phi_\sigma(x)=\Phi(x).
\end{equation}
By definition of the imploded moment map, the composition of the
bottom two arrows in \eqref{equation;moment-commute} is the map
sending $x\in\Phi\inv(\exp\sigma)$ to $\Phi(x)$.  Precomposing with
the covering map
$\tilde{\Phi}_\sigma\inv(K_2)\to\Phi\inv(\exp\sigma)$, which sends
$(x,k_1,k_2)$ to $x$, we get the map $(x,1,k_2)\mapsto\Phi(x)$, which
is the same as \eqref{equation;arrows}.
\end{proof}

Let $\{\,X_i\mid i\in I\,\}$ be the collection of all strata of all
the pieces $X_\sigma$, where $\sigma$ ranges over all faces of the
alcove.  The imploded cross-section is their disjoint union,
\begin{equation}\label{equation;stratification}
M\impl=\coprod_{i\in I}X_i,
\end{equation}
and we shall call the $X_i$ the \emph{strata} of $M\impl$ (although we
shall not prove here that they form a stratification of $M\impl$ in
the technical sense).  The preceding results can be summarized as
follows.

\begin{theorem}\label{theorem;decomposition}
The decomposition \eqref{equation;stratification} of the imploded
cross-section is a locally finite partition into locally closed
subspaces, each of which is a symplectic manifold.  There is a unique
open stratum, which is dense in $M\impl$ and symplectomorphic to the
principal cross-section of $M$.  The action of the maximal torus $T$
on $M\impl$ preserves the decomposition and the imploded moment map
$\Phi\impl\colon M\impl\to T$ restricts to a moment map for the
$T$-action on each stratum.
\end{theorem}

For this reason we call $M\impl$ a \emph{stratified quasi-Hamiltonian
$T$-space}.  Quotients of $M\impl$ by the $T$-action can be defined
just as if $M\impl$ was a manifold: for
$g\in\Phi(M)\cap\exp\bar{\alcove}=\Phi\impl(M\impl)$ one puts
$M\impl\quot[g]T=\Phi\impl\inv(g)/T$.  Again these quotients decompose
into symplectic manifolds according to orbit type.  This is seen by
applying Theorem \ref{theorem;singular-quotient} to all of the
$T$-manifolds $X_i$ that intersect $\Phi\impl\inv(g)$.  Moreover,
$\Phi\impl\inv(g)$ is the quotient of $\Phi\inv(g)$ by $[K_g,K_g]$ and
the quotient map $\Phi\inv(g)\to\Phi\impl\inv(g)$ descends to a
homeomorphism $M\quot[g]K\to M\impl\quot[g]T$.  It is now
straightforward to check the following assertion.

\begin{addendum}\label{addendum;quotient}
For all $g\in\Phi(M)\cap\exp\bar{\alcove}$, the homeomorphism
$M\quot[g]K\to M\impl\quot[g]T$ induced by the $[K_g,K_g]$-orbit map
$\Phi\inv(g)\to\Phi\impl\inv(g)$ sends strata to strata and restricts
to a symplectomorphism on each stratum.
\end{addendum}

Thus the imploded cross-section $M\impl$ is the ``abelianization'' of
$M$ in the sense that every quasi-symplectic quotient of $M$ by the
$K$-action can be written as a quasi-symplectic quotient of $M\impl$
by the $T$-action.

Just like quasi-Hamiltonian reduction, quasi-Hamiltonian implosion for
the action of a product $K_1\times K_2$ can be performed with respect
to one of the factors, say $K_2$.  Using Remark \ref{remark;product}
one can show the resulting space is quasi-Hamiltonian for the product
$K_1\times T_2$, where $T_2$ is a maximal torus of $K_2$.  We shall
not develop this variant of the theory save in the important special
case of the double $DK$, which is investigated in the next section.

\section{The universal imploded cross-section}\label{section;double}

\subsection*{Implosion of the double}

The cotangent bundle $T^*K$ is a Hamiltonian $K\times K$-manifold.
Implosion with respect to, say, the right $K$-action yields a $K\times
T$-space $T^*K\impl$, which stratifies into symplectic manifolds.
This space is the \emph{universal imploded cross-section} in the sense
that the imploded cross-section of every Hamiltonian $K$-manifold $M$
is a quotient,
$$
M\impl\cong(M\times T^*K\impl)\quot K.
$$
Here the quotient is taken with respect to the given $K$-action on $M$
and the residual (left) $K$-action on $T^*K\impl$.  (See \cite[Theorem
4.9]{guillemin-jeffrey-sjamaar}.)  A similar phenomenon occurs in the
quasi-Hamiltonian category.  According to \cite[Section
3.2]{alekseev-malkin-meinrenken;lie-group} the quasi-Hamiltonian
analogue of the cotangent bundle $T^*K$ is the double $DK=K\times
K$\glossary{D@$DK$, double $K\times K$}.  The $K\times K$-action on
$DK$ is defined by $(g_1, g_2)(u,v)=(g_1ug_2\inv,\Ad(g_2)v)$, the
two-form by
\begin{equation}\label{equation;form}
\omega=-\frac12\bigl(\Ad(v)u^*\theta_L,u^*\theta_L\bigr)
-\frac12\bigl(u^*\theta_L,v^*(\theta_L+\theta_R)\bigr),
\end{equation}
and the moment map $\Psi\colon DK\to K\times K$ by
$\Psi=\Psi_1\times\Psi_2$ with
\begin{equation}\label{equation;moment-map}
\Psi_1(u,v)=\Ad(u)v\inv\qquad\text{and}\qquad\Psi_2(u,v)=v.
\end{equation}
(Here we work in the coordinate system of \cite[Remark
3.2]{alekseev-malkin-meinrenken;lie-group}.  Actually
[\emph{loc.\ cit.}] uses the opposite space $DK^-=(K\times
K,-\omega,i\circ\Psi)$, but $DK^-$ is isomorphic to $DK$ via the
diffeomorphism $(u,v)\mapsto(u,v\inv)$.)


The fusion product $M\fusion DK$ (see Section
\ref{section;hamiltonian}) of an arbitrary quasi-Hamiltonian
$K$-manifold $M$ and the quasi-Hamiltonian $K\times K$-manifold $DK$
is a quasi-Hamiltonian $K\times K$-manifold.  Likewise, for each face
$\sigma\le\alcove$ the cross-section
$\Psi_2\inv(\eu{S}_\sigma)=K\times\eu{S}_\sigma$ of $DK$ is a
quasi-Hamiltonian $K\times K_\sigma$-manifold (see Remark
\ref{remark;product}) and so $M\fusion(K\times\eu{S}_\sigma)$ is a
quasi-Hamiltonian $K\times K_\sigma$-manifold.

Now define $j\colon M\to M\fusion DK$ by $j(m)=(m,1,\Phi(m))$.  The
following result is analogous to \cite[Lemma
4.8]{guillemin-jeffrey-sjamaar} and is proved in the same way.  (An
embedding is \emph{quasi-symplectic} if it preserves the two-forms.)

\begin{lemma}\label{lemma;universal}
\begin{enumerate}
\item\label{item;universal-centre}
The map $j$ is a quasi-symplectic embedding and induces an isomorphism
of quasi-Hamiltonian $K$-manifolds
$$
\bar\jmath\colon M\to(M\fusion DK)\quot K.
$$
Here the right-hand side is the quotient with respect to the diagonal
$K$-action, with $K$ acting on the left on $DK$.  The $K$-action on
$(M\fusion DK)\quot K$ is the one induced by the right action on $DK$.
\item\label{item;universal}
For every $\sigma\le\alcove$, $j$ maps $M_\sigma$ into
$M\fusion(K\times\eu S_\sigma)$ and induces an isomorphism of
quasi-Hamiltonian $K_\sigma$-manifolds
$$
\bar\jmath_\sigma\colon M_\sigma\to\bigl(M\fusion(K\times\eu
S_\sigma)\bigr)\bigquot K,
$$
where the quotient is taken as in \eqref{item;universal-centre}.
\end{enumerate}
\end{lemma}

We now implode $DK$, but using only the right action of $K$ and the
second component $\Psi_2$ of the moment map.  In other words we form
the quotient topological space of
$\Psi_2\inv(\exp\bar{\alcove})=K\times\exp\bar{\alcove}$ by the
equivalence relation $(u,v)\sim(ug\inv,v)$ for $g\in[K_v,K_v]$.  As in
\eqref{equation;quasi-decomposition} the resulting space is a union of
subspaces,
\begin{equation}\label{equation;universal-decomposition}
DK\impl
=\coprod_{\sigma\le\alcove}K/[K_\sigma,K_\sigma]\times\exp\sigma.
\end{equation}
The moment map $\Psi_2$ is transverse to all faces of the alcove, so
it happens that each piece
$X_\sigma=K/[K_\sigma,K_\sigma]\times\exp\sigma$ in this partition is
a smooth manifold.  We shall show that $X_\sigma$ is a
quasi-Hamiltonian $K\times T$-manifold by writing it as a
quasi-symplectic quotient up to a covering.  Namely let
$DK_\sigma=\Psi_2\inv(\eu{S}_\sigma)=K\times\eu{S}_\sigma$.  By Remark
\ref{remark;product} $DK_\sigma$ is a Hamiltonian $K\times
K_\sigma$-manifold with moment map $\Psi_\sigma$ being the restriction
of $\Psi$ to $DK_\sigma$.  Let $p\colon\tilde{K}_\sigma\to K_\sigma$
be the covering \eqref{equation;covering} of the centralizer
$K_\sigma$.  The pullback of this covering via $\Psi_\sigma$ is
$DK\sptilde_\sigma=K\times\tilde{\eu{S}}_\sigma$, where
$\tilde{\eu{S}}_\sigma=p\inv(\eu{S}_\sigma)$.  Thus
$K\times\tilde{\eu{S}}_\sigma$ is a quasi-Hamiltonian
$K\times\tilde{K}_\sigma$-manifold.  Its quotient with respect to the
$[K_\sigma,K_\sigma]$-action is the quasi-Hamiltonian $K\times
Z(K_\sigma)^0$-manifold
$$
\tilde{X}_\sigma=(K\times\tilde{\eu{S}}_\sigma)
\quot{[K_\sigma,K_\sigma]} =K/[K_\sigma,K_\sigma]\times
p\inv(\exp\sigma),
$$
which is a covering space of $X_\sigma$ with covering group
$\Gamma_\sigma$.  The following is an analogue of Corollary
\ref{corollary;stratification} and Lemma \ref{lemma;implosion} and is
proved in the same way.

\begin{lemma}\label{lemma;universal-cover}
For every $\sigma\le\alcove$ the subspace
$X_\sigma=K/[K_\sigma,K_\sigma]\times\exp\sigma$ of $DK\impl$ is a
quasi-Hamiltonian $K\times T$-manifold.  The moment map $X_\sigma\to
K\times T$ is the restriction to $X_\sigma$ of the continuous map
$\Psi\impl\colon DK\impl\to K\times T$ induced by $\Psi\colon DK\to
K\times K$.
\end{lemma}

We can now establish the universal property of $DK\impl$.  Consider an
arbitrary quasi-Hamiltonian $K$-manifold $(M,\omega,\Phi)$.  We form
the fusion product $M\fusion DK\impl$ just as if $DK\impl$ was a
smooth manifold: the underlying space is the product $M\times
DK\impl$, the action and the moment map are defined as in the smooth
case, and the two-form is defined stratum by stratum.  The resulting
space is a stratified space whose strata are quasi-Hamiltonian
$K\times T$-manifolds.  Stratum by stratum we can reduce $M\fusion
DK\impl$ with respect to the $K$-action to obtain a stratified
quasi-Hamiltonian $T$-space $(M\fusion DK\impl)\quot K$.  This
quotient turns out to be the same as $M\impl$.

\begin{theorem}[universality of the imploded double]%
\label{theorem;universal}
Let $M$ be a quasi-Hamiltonian $K$-manifold.  The map $j$ induces a
homeomorphism
$$
j\impl\colon M\impl\overset{\cong}{\rightarrow}(M\fusion DK\impl)\quot
K,
$$
which maps strata to strata and whose restriction to each stratum is
an isomorphism of quasi-Hamiltonian $T$-manifolds.
\end{theorem}

\begin{proof}
Observe that $j$ maps $\Phi\inv(\exp\alcove)$ into
$M\times\Phi_2\inv(\exp\alcove)$ and that $m_1\sim m_2$ implies
$j(m_1)=j(m_2)$ for all $m_1$, $m_2\in\Phi\inv(\exp\alcove)$.  Hence
$j$ descends to a continuous map $j'\colon M\impl\to M\fusion
DK\impl$.  The $K\times T$-moment map on $M\fusion DK\impl$ is
$\phi=\phi_1\times\phi_2$, where
$\phi_1\bigl(m,(u,v)\mod{\sim}\bigr)=\Phi(m)uv\inv u\inv$ and
$\phi_2\bigl(m,(u,v)\mod{\sim}\bigr)=v$.  Therefore $j'$ maps $M\impl$
into $\phi_1\inv(1)$.  Dividing by the action of $K$ we obtain the map
$j\impl$.  It is straightforward to check that $j\impl$ is bijective
and therefore a homeomorphism.  Moreover $j\impl$ maps the piece
$\Phi\inv(\exp\sigma)/[K_\sigma,K_\sigma]$ of $M\impl$ onto $(M\fusion
X_\sigma)\quot K$ (where $X_\sigma$ is as in Lemma
\ref{lemma;universal-cover}).  It remains to show that the restriction
of $j\impl$,
\begin{equation}\label{equation;isomorphism}
\Phi\inv(\exp\sigma)/[K_\sigma,K_\sigma]\to(M\fusion X_\sigma)\quot K,
\end{equation}
maps strata to strata and is an isomorphism of quasi-Hamiltonian
$T$-manifolds on each stratum.  To see this, recall the isomorphism of
Lemma \ref{lemma;universal}\eqref{item;universal}, which lifts to an
isomorphism of quasi-Hamiltonian $\tilde{K}_\sigma$-manifolds
$$
\tilde{M}_\sigma\to\bigl(M\fusion(K\times\tilde{\eu
S}_\sigma)\bigr)\bigquot K.
$$
Upon reduction by $[K_\sigma,K_\sigma]$ this yields a
stratification-preserving homeomorphism
$$
\tilde{M}_\sigma\quot{[K_\sigma,K_\sigma]}
\to(M\fusion\tilde{X}_\sigma)\quot K,
$$
whose restriction to each stratum is an isomorphism of $K\times
Z(K_\sigma)^0$-spaces.  Taking the quotient by $\Gamma_\sigma$ on both
sides we obtain the map \eqref{equation;isomorphism}.  It now follows
from Corollary \ref{corollary;stratification} that this map sends
strata isomorphically onto strata.
\end{proof}

\subsection*{Discrete symmetries}

The centre of $K$ gives rise to a set of discrete symmetries of the
universal imploded cross-section $DK\impl$.  A stratum
$K/[K_\sigma,K_\sigma]\times\exp\sigma$ of $DK\impl$ consists of a
single point if and only if $\sigma$ is a vertex of the alcove and
$\exp\sigma$ is an element of the centre $Z(K)$.  Thus one-point
strata are in one-to-one correspondence with the centre.  These strata
are exactly the fixed points under the $K\times T$-action.  It turns
out there exists an extension $\hat{T}$ of the centre by the maximal
torus that acts on $DK\impl$, permutes the one-point strata, and in
addition preserves the quasi-Hamiltonian structure.  This action has
no Hamiltonian analogue, because $T^*K\impl$ has just a single
one-point stratum (the ``vertex'').  The material below relies on
facts concerning the action of the centre on the alcove, which are
reviewed at the end of Appendix \ref{section;adjoint}.

Let $\zeta\colon Z(K)\to W$ be the injective homomorphism of Lemma
\ref{lemma;centre} and let $\hat{T}$\glossary{T@$\hat{T}$, extension
of centre by maximal torus} be the preimage of $\zeta(Z(K))$ under the
canonical projection $N(T)\to W$, where $N(T)$ is the normalizer of
the maximal torus.  Thus $\hat{T}$ is a subgroup of $N(T)$ and we have
short exact sequence
\begin{equation}\label{equation;exact-centre}
\xymatrix{{\smash{T}}\vphantom{Z(K)}\ar@{((->}[r]&{\smash{\hat{T}}}
\vphantom{Z(K)}\ar@{->>}[r]^-\pi&Z(K).}
\end{equation}

\begin{example}\label{example;exact-centre}
The centre of $K=\SU(l+1)$ is generated by $c=e^{-2\pi
i/(l+1)}I_{l+1}$, where $I_{l+1}$ is the identity matrix, and by
Example \ref{example;centre-weyl} $\zeta(c)$ is the cyclic permutation
$x\mapsto(x_{l+1},x_1,x_2\dots,x_l)$.  A representative of $c$ in
$\hat{T}$ is the matrix
$$\hat{c}=-\begin{pmatrix}0&(-1)^l\\I_l&0\end{pmatrix},$$
%
and therefore $\hat{T}=\bigcup_{j=0}^l\hat{c}^jT$.  Note that
$\hat{T}$ is not abelian.  For $t\in T$ the order of $\hat{c}t$ is
$l+1$ if $l$ is even and $2l+2$ if $l$ is odd.  Hence the sequence
\eqref{equation;exact-centre} splits if and only if $l$ is even.
\end{example}

Define a smooth $\hat{T}$-action $\eu{A}$\glossary{A@$\eu{A}$,
$\hat{T}$-action on $DK$} on $DK$ by
$$
\eu{A}_g(u,v)=\bigl(ug\inv,\pi(g)\Ad(g)v\bigr)
=\bigl(ug\inv,g\pi(g)vg\inv\bigr).
$$
Since $T=\ker\pi$ this action is an extension of the right $T$-action
on $DK$.

\begin{lemma}\label{lemma;centre-action}
The $\hat{T}$-action preserves
$\Psi_2\inv(\exp\bar{\alcove})=K\times\exp\bar{\alcove}$ and induces a
continuous effective action on $DK\impl$.  This induced action
commutes with the $K$-action, maps strata to strata, and is transitive
on the collection of one-point strata.
\end{lemma}

\begin{proof}
If $v\in\exp\bar{\alcove}$ then $\pi(g)gvg\inv\in\exp\bar{\alcove}$ by
Lemma \ref{lemma;centre} (applied to $w=gT\in W$ and $c=\pi(g)$), so
$\eu{A}_g$ maps $\Psi_2\inv(\exp\bar{\alcove})$ to itself.  Let
$(u,v)\in K\times\exp\bar{\alcove}$ and $k\in[K_v,K_v]$.  Then $(u,v)$
is equivalent to $(uk\inv,v)$ and
$$
\eu{A}_g(uk\inv,v)=\bigl(uk\inv g\inv,\pi(g)\Ad(g)v\bigr)
=\bigl(ug\inv\Ad(g)k\inv,\pi(g)\Ad(g)v\bigr),
$$
where
$$
\Ad(g)k\in\Ad(g)[K_v,K_v]
=\bigl[K_{\pi(g)\Ad(g)v},K_{\pi(g)\Ad(g)v}\bigr],
$$
so $\eu{A}_g(u,v)\sim f_g(uk\inv,v)$.  Hence $\eu{A}_g$ descends to a
homeomorphism $\bar{\eu{A}}_g$ from $DK\impl$ to itself.  On the open
stratum $K\times\exp\bar{\alcove}$ $g$ acts by right translation on
$K$ and via the central element $\pi(g)$ on $\alcove$.  This shows
that the induced action is effective.  Since the action $\eu{A}$
commutes with the left $K$-action on $DK$, $\bar{\eu{A}}$ commutes
with the induced $K$-action on $DK\impl$.  The action of the coweight
lattice $P(R\spcheck)=\exp_T\inv(Z(K))$ on $\t$ is affine and maps
alcoves to alcoves.  (See the proof of Lemma \ref{lemma;centre}.)
Therefore it maps faces of alcoves to faces of alcoves.  This implies
that $\bar{\eu{A}}_g$ maps strata to strata.  A one-point stratum is
an equivalence class $(1,c)\mod{\sim}$ for some $c\in Z(K)$.  Since
$\eu{A}_g(1,c)=(g\inv,\pi(g)c)\sim(1,\pi(g)c)$, we have
$\bar{\eu{A}}_g\bigl((1,c)\mod{\sim}\bigr)
=\bigl((1,\pi(g)c)\mod{\sim}\bigr)$, so $\hat{T}$ acts transitively on
the set of one-point strata.
\end{proof}

Observe however that the centre $Z(K)$ itself does not act on
$DK\impl$ (unless the sequence \eqref{equation;exact-centre} splits).
Also, $DK\impl$ is not quite a quasi-Hamiltonian $K\times
\hat{T}$-manifold.  The conjugation action of $K\times\hat{T}$ on
$K\times T$ is given by
$$\Ad(k,g)\cdot(g_1,g_2)=\bigl(\Ad(k)g_1,\Ad(g)g_2\bigr),$$
whereas the imploded moment map $\Psi\impl$ is equivariant only under
the ``affine'' conjugation action $\func{Ad\sphat\,}(k,g)(g_1,g_2)
=\bigl(\pi(g)\inv\Ad(k)g_1,\pi(g)\Ad(g)v\bigr)$.

\begin{proposition}\label{proposition;centre-action}
The $K\times\hat{T}$-action on $DK\impl$ preserves the two-forms on
the strata.  The $K\times T$-moment map $\Psi\impl\colon DK\impl\to
K\times T$ is $\func{Ad\sphat\,}$-equivariant for the
$K\times\hat{T}$-action.
\end{proposition}

\begin{proof}
The map $\eu{A}_g$ is the composition of
$(u,v)\mapsto(ug\inv,\Ad(g)v)$ and $(u,v)\mapsto(u,\pi(g)v)$.  The
first map is the right action of the element $g$ and therefore
preserves $\omega$; the second map preserves $\omega$ because both
$\theta_L$ and $\theta_R$ are invariant under translation by the
central element $\pi(g)$.  It follows that $\bar{\eu{A}}_g$ preserves
the forms on the strata.  For $(k,g)\in K\times\hat{T}$ and $(u,v)\in
K\times\exp\bar{\alcove}$ we have
\begin{multline*}
\Psi\bigl((k,g)\cdot(u,v)\bigr)=\Psi\bigl(kug\inv,\pi(g)\Ad(g)v\bigr)\\
=\bigl(\pi(g)\inv \Ad(ku)v\inv,\pi(g)\Ad(g)v\bigr)
=\func{Ad\sphat\,}(k,g)\Psi(u,v),
\end{multline*}
which establishes the equivariance of $\Psi$.
\end{proof}

Further discrete symmetries of the universal imploded cross-section
come from the duality automorphism of the alcove, which is defined by
$\xi\mapsto-w_0\xi$, where $w_0$\glossary{w@$w_0$, longest Weyl group
element} is the longest Weyl group element.  Let
$T_0$\glossary{T@$T_0$, extension of $\Z/2\Z$ by maximal torus} be the
inverse image of the subgroup $\langle w_0\rangle=\{1,w_0\}$ of $W$
under the projection $N(T)\to T$.  It fits into an exact sequence
\begin{equation}\label{equation;exact-long}
\xymatrix{{\smash{T}}\vphantom{Z(K)}\ar@{((->}[r]&{\smash{T_0}}
\vphantom{Z(K)}\ar@{->>}[r]&{\smash{\langle
w_0\rangle\cong\Z/2\Z}}\vphantom{Z(K)}.}
\end{equation}
Define a smooth $T_0$-action on $DK$ by $\eu{B}_g(u,v)
=\bigl(ug\inv,\Ad(g)v^{\rho(g)}\bigr)$,\glossary{B@$\eu{B}$,
$T_0$-action on $DK$} where $\rho(g)=1$ if $g$ projects to $1\in W$
and $\rho(g)=-1$ if $g$ projects to $w_0$.  The following result is
analogous to Lemma \ref{lemma;centre-action} and Proposition
\ref{proposition;centre-action}.

\begin{proposition}\label{proposition;duality}
The $T_0$-action preserves
$\Psi_2\inv(\exp\bar{\alcove})=K\times\exp\bar{\alcove}$ and descends
to a continuous effective action on $DK\impl$.  This action commutes
with the $K$-action and maps strata to strata.  An element $g\in T_0$
preserves or reverses the two-forms on the strata according to whether
$\rho(g)=1$ or $\rho(g)=-1$.  The $K\times T$-moment map
$\Psi\impl\colon DK\impl\to K\times T$ is $T_0$-equivariant for the
$T_0$-action $g\cdot(g_1,g_2)
=\bigl(g_1^{\rho(g)},\Ad(g)g_2^{\rho(g)}\bigr)$ on $K\times T$.
\end{proposition}

\begin{example}\label{example;exact-long}
The longest Weyl group element of $K=\SU(l+1)$ is the permutation
$w_0(x)=(x_{l+1},x_l,\dots,x_1)$.  One checks readily that $w_0$ has a
representative of order two in $N(T)$ if and only if
$l\not\equiv1\mod4$; if $l\equiv1\mod4$ it has a representative of
order four.  Therefore the sequence \eqref{equation;exact-long} splits
if and only if $l\not\equiv1\mod4$.  A convenient choice for such a
representative is
$$
n_0=
\begin{pmatrix}0&(-1)^{l(l+1)/2}\\J_l&0\end{pmatrix}
\quad\text{if $l\equiv0$, $1$, $3\mod4$,}
$$
where $J_l$ is the antidiagonal $l\times l$-matrix with ones up the
second diagonal, and
$$
n_0=
\begin{pmatrix}0&0&J_{2m+1}\\0&-1&0\\J_{2m+1}&0&0\end{pmatrix}
\quad\text{if $l=4m+2\equiv2\mod4$.}
$$
For $l=1$ we have $T_0=\hat{T}=N(T)$, but the actions of $T_0$ and
$\hat{T}$ on $DK\impl$ are not the same.  
\end{example}

\subsection*{A smoothness criterion}

Certain strata of the universal imploded cross-section are redundant
in the sense that $DK\impl$ is a smooth manifold at these points and
the two-form and the moment map extend smoothly.  Thus these
singularities are removable in the strongest possible sense.  For
instance if $K$ is a product of copies of $\SU(2)$ \emph{all}
singularities are removable.  To prove this we must be more explicit
about the analogy between the double $DK$ and the cotangent bundle
$T^*K$.

Identify $\k^*$ with $\k$ by means of the inner product, trivialize
$T^*K$ by means of left translations, and let $\ca{H}$ be the $K\times
K$-equivariant map $\id\times\exp\colon T^*K\cong K\times\k\to DK$.
Let $\omega_0$ be the standard symplectic form on $T^*K$ and
$\Psi_0(g,\lambda)=(-\Ad(g)\lambda,\lambda)$ the standard moment map
for the $K\times K$-action.  Let $O$ be the set of all $\xi\in\t$ with
$(2\pi i)\inv\alpha(\xi)<1$ for all roots $\alpha$ and let
$U=\Ad(K)O$.

\begin{proposition}\label{proposition;cotangent-double}
The triple $(T^*K,\ca{H}^*\omega,\ca{H}^*\Psi)$ is the exponentiation
(see Section \ref{section;hamiltonian}) of the Hamiltonian $K\times
K$-manifold $(T^*K,\omega_0,\Psi_0)$.  Hence $(K\times
U,\ca{H}^*\omega,\ca{H}^*\Psi)$ is a quasi-Hamiltonian $K\times
K$-manifold.
\end{proposition}

\begin{proof}
A tangent vector to $T^*K=K\times\k$ at $(g,\lambda)$ is of the form
$(L(g)_*\xi,\eta)$ with $\xi$, $\eta\in\k$.  The standard symplectic
form is
\begin{equation}\label{equation;omega0}
(\omega_0)_{(g,\lambda)}
\bigl((L(g)_*\xi_1,\eta_1),(L(g)_*\xi_2,\eta_2)\bigr)
=(\xi_2,\eta_1)-(\xi_1,\eta_2).
\end{equation}
A calculation using \eqref{equation;maurer} and \eqref{equation;form}
yields
\begin{multline}\label{equation;omega}
\ca{H}^*\omega_{(g,\lambda)}
\bigl((L(g)_*\xi_1,\eta_1),(L(g)_*\xi_2,\eta_2)\bigr)
=-\bigl((\sinh\ad\lambda)\xi_1,\xi_2\bigr)\\
-\biggl(\frac{\sinh\ad\lambda}{\ad\lambda}\xi_1,\eta_2\biggr)
+\biggl(\frac{\sinh\ad\lambda}{\ad\lambda}\xi_2,\eta_1\biggr).
\end{multline}
Let $\pi_1$, $\pi_2\colon\k\oplus\k\to\k$ be the projections onto the
respective factors and let $\varpi'=\pi_1^*\varpi+\pi_2^*\varpi$,
where $\varpi\in\Omega^2(\k)^K$ is given by \eqref{equation;varpi}.
From
$$
T_{(g,\lambda)}\Psi_0(L(g)_*\xi,\eta)
=\bigl(-\Ad(g)([\xi,\lambda]+\eta),\eta\bigr)
$$
we obtain
\begin{multline}\label{equation;pull-varpi}
\Psi_0^*\varpi'_{(g,\lambda)}
\bigl((L(g)_*\xi_1,\eta_1),(L(g)_*\xi_2,\eta_2)\bigr)\\
=\varpi_{-\lambda}([\xi_1,\lambda]+\eta_1,[\xi_2,\lambda]+\eta_2)
+\varpi_\lambda(\eta_1,\eta_2)\\
\,\quad=\bigl((\ad\lambda-\sinh\ad\lambda)\xi_1,\xi_2\bigr)
+(\xi_1,\eta_2)-(\xi_2,\eta_1)\\
-\biggl(\frac{\sinh\ad\lambda}{\ad\lambda}\xi_1,\eta_2\biggr)
+\biggl(\frac{\sinh\ad\lambda}{\ad\lambda}\xi_2,\eta_1\biggr).
\end{multline}
From \eqref{equation;omega0}--\eqref{equation;pull-varpi} we conclude
$\ca{H}^*\omega=\omega_0+\Psi_0^*\varpi'$.  Let $\exp\sptilde$ denote
the exponential map of $K\times K$.  Then
$\exp\sptilde(\Psi_0(g,\lambda))
=\bigl(\Ad(g)\exp(-\lambda),\exp\lambda\bigr)$ and therefore
$\ca{H}^*\Psi=\exp\sptilde\circ\Psi_0$.  Lastly by Lemma
\ref{lemma;chart}\eqref{item;etale-chart} $\exp\colon U\to K$ is a
local diffeomorphism onto its image.  Therefore $(K\times
U,\ca{H}^*\omega,\ca{H}^*\Psi)$ is a quasi-Hamiltonian $K\times
K$-manifold.
\end{proof}

The local diffeomorphism $\ca{H}=\id\times\exp\colon K\times U\to
K\times\exp U$ descends to the imploded cross-sections as follows.
Under the isomorphism $\k^*\cong\k$ the chamber $\ca C$ is mapped to
the dual chamber $\ca{C}\spcheck$.  By Lemma
\ref{lemma;chart}\eqref{item;chamber-alcove} the intersection
$\bar{\ca C}\cap O$ is equal to $\bar\alcove\cap O$ and therefore
$\ca{H}$ maps $\Psi_0\inv(\bar{\ca C})\cap(K\times U)$ onto
$\Psi\inv(\bar\alcove)\cap(K\times\exp U)$.  Moreover
$K_\xi=K_{\exp\xi}$ for $\xi\in\bar\alcove\cap O$ by Lemma
\ref{lemma;face}\eqref{item;face-face}, so $\ca{H}$ induces a $K\times
T$-equivariant surjective local homeomorphism
\begin{equation}\label{equation;homeo}
h\colon(K\times U)\impl\to(K\times\exp U)\impl.
\end{equation}
The left-hand side is open in the Hamiltonian implosion $T^*K\impl$
and the right-hand side is an open and dense subset of the
quasi-Hamiltonian implosion $DK\impl$, namely the union of the strata
corresponding to the faces $\sigma\le\alcove$ satisfying $\sigma\ge0$.
In fact $h$ is injective and therefore a homeomorphism.  To see
this, let $\sigma\ge0$ and let $\sigma'\le\ca{C}$ be the unique face
of the chamber containing $\sigma$.  Then $\sigma'\cap O=\sigma$ and
$K_{\sigma'}=K_\sigma$ and the restriction of $h$ to the stratum
corresponding to $\sigma'$ is the bijection
$$
h=\id\times\exp\colon K/[K_\sigma,K_\sigma]\times\sigma\to
K/[K_\sigma,K_\sigma]\times\exp\sigma.
$$
We can now show that certain rather special singularities of $DK\impl$
are removable.

\begin{theorem}\label{theorem;smooth}
Let $\sigma$ be a face of $\alcove$ satisfying
$[K_\sigma,K_\sigma]\cong\SU(2)^k$ (resp.\
$[\k_\sigma,\k_\sigma]\cong\lie{su}(2)^k$) for some $k$ and possessing
a vertex $\xi$ such that $\exp\xi$ is central.  Then $DK\impl$ is a
smooth quasi-Hamiltonian $K\times T$-manifold (resp.\ orbifold) in a
neighbourhood of the stratum corresponding to $\sigma$.
\end{theorem}

\begin{proof}
After acting by the central element $\exp-\xi$ we may assume $\xi=0$
and $\sigma\ge0$.  Then the stratum of $DK\impl$ associated with
$\sigma$ is contained in the image $(K\times\exp U)\impl$ of $h$.  Let
$\sigma'$ be the face of the chamber containing $\sigma$.  Put
$X=DK\impl$ and $X'=T^*K\impl$.  Let
$X_\sigma=K/[K_\sigma,K_\sigma]\times\exp\sigma$ be the stratum of $X$
associated with $\sigma$ and let
$X'_\sigma=K/[K_\sigma,K_\sigma]\times\sigma'$ be the stratum of $X'$
associated with $\sigma'$.  Let $(\omega_\sigma,\Psi_\sigma)$ denote
the quasi-Hamiltonian structure on $X_\sigma$ and
$(\omega'_\sigma,\Psi'_\sigma)$ the Hamiltonian structure on
$X'_\sigma$.  It follows from Proposition
\ref{proposition;cotangent-double} that
\begin{equation}\label{equation;linearization-stratum}
h^*\Psi_\sigma=\exp\sptilde\circ\Psi'_\sigma
\qquad\text{and}\qquad
\bar f^*\omega_\sigma=\omega'_\sigma+(\Psi'_\sigma)^*\varpi,
\end{equation}
where $\exp\sptilde$ is now the exponential map for $K\times T$.  By
\cite[Proposition 6.15]{guillemin-jeffrey-sjamaar} the singularities
of $T^*K\impl$ are removable at $X'_\sigma$ in the sense that the open
set $Y=\bigcup_{\tau\ge\sigma}X'_\tau$ possesses a smooth structure
(resp.\ an orbifold structure if
$[\k_\sigma,\k_\sigma]\cong\lie{su}(2)^k$) and the Hamiltonian
structure on the open stratum extends smoothly to this open set.  With
respect to this smooth structure the strata $X'_\tau$ for
$\tau\ge\sigma$ are smooth submanifolds of $Y$ and the Hamiltonian
structure on $Y$ restricts to the given Hamiltonian structures on the
strata.  From \eqref{equation;linearization-stratum} we then see that
the quasi-Hamiltonian structure likewise extends smoothly to $Y$.  Via
the homeomorphism $h$ we now obtain a smooth quasi-Hamiltonian
$K\times T$-manifold (resp.\ orbifold) structure on
$h(Y)=\bigcup_{\tau\ge\sigma}X_\tau$.
\end{proof}

\subsection*{Examples}

Next we describe in more detail the imploded cross-sections of the
doubles of $\SU(n)$ and $\SU(2,\H)$.  It turns out that certain strata
in these spaces have a nonsingular closure.  This rather serendipitous
observation brings to light a family of quasi-Ham\-ilton\-ian
manifolds, which surprisingly turns out to be new.  (See Theorem
\ref{theorem;s2n}.)  In each of the examples below we specify an
ordering $\sigma_0$, $\sigma_1$,\dots, $\sigma_r$ of the vertices of
the alcove.  The face corresponding to a subset $I$ of
$\{0,1,\dots,r\}$ is then denoted by $\sigma_I$, the centralizer of
$\sigma_I$ by $K_I$ and the associated stratum of $X=DK\impl$ by
$X_I$.

\begin{example}\label{example;su3}
Let $K=\SU(3)$ and let $T$ be its standard maximal torus.  The alcove
is an open equilateral triangle (cf.\ Example
\ref{example;centre-weyl}) with three vertices $\sigma_j$,
corresponding to the elements $c^j$ of the centre of $K$.  The
one-dimensional faces are the open segments $\sigma_{jk}$ joining
$\sigma_j$ to $\sigma_k$, where $0\le j<k\le2$.
$$
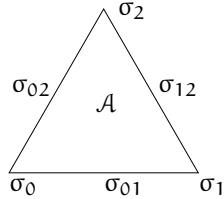
\setlength{\unitlength}{0.003mm}
\begin{picture}(9046,8013)(0,-10)
\path(8400,654)(4211,7910)(22,654)(8400,654)
\put(0,-150){\makebox(0,0)[lb]{\smash{$\sigma_0$}}}
\put(8400,-150){\makebox(0,0)[lb]{\smash{$\sigma_1$}}}
\put(4800,7654){\makebox(0,0)[lb]{\smash{$\sigma_2$}}}
\put(3900,3354){\makebox(0,0)[lb]{\smash{$\alcove$}}}
\put(4200,-150){\makebox(0,0)[lb]{\smash{$\sigma_{01}$}}}
\put(6600,4254){\makebox(0,0)[lb]{\smash{$\sigma_{12}$}}}
\put(100,4254){\makebox(0,0)[lb]{\smash{$\sigma_{02}$}}}
\end{picture}
$$
The centralizers and their commutator subgroups are
\begin{alignat*}{3}
&K_\alcove=T,&\qquad&K_{jk}\cong\group{S}(\U(2)\times\U(1)),
&\qquad&K_j=K,\\
&[K_\alcove,K_\alcove]=\{I\},&\qquad&[K_{jk},K_{jk}]\cong\SU(2),
&\qquad&[K_j,K_j]=K.
\end{alignat*}
Thus the open stratum $X_\alcove$ is diffeomorphic to
$\SU(3)\times\alcove$ and has real dimension $10$; the intermediate
strata $X_{jk}$ are diffeomorphic to $\SU(3)/\SU(2)\times(0,1)\cong
S^5\times(0,1)$; and the smallest strata $X_j$ are points in
one-to-one correspondence with the centre.  By Theorem
\ref{theorem;smooth} the $X_{jk}$ are removable singularities.  By
\cite[Example 6.16]{guillemin-jeffrey-sjamaar} $T^*K\impl$ is
isomorphic to the quadric cone
$$\{\,(y,x)\in(\C^3)^*\times \C^3\mid y(x)=0\,\}$$
with a K\"ahler structure obtained by restricting the ambient flat
K\"ahler structure on $(\C^3)^*\times \C^3$.  Because of the local
isomorphism between $T^*K\impl$ and $DK\impl$ we conclude that the
three smallest strata of $DK\impl$ are genuine singularities.  The
intermediate strata in $T^*K\impl$ are the punctured axes,
$\{\,(y,0)\in(\C^3)^*\times \C^3\mid y\ne0\,\}$ and
$\{\,(0,x)\in(\C^3)^*\times \C^3\mid x\ne0\,\}$.  The closure of each
is a copy of a flat $\C^3$.  By the argument given in the proof of
Theorem \ref{theorem;smooth} we conclude that the closure of each
$X_{jk}$ is nonsingular, namely a six-sphere.  The upshot is that
$S^6$ is a quasi-Hamiltonian $\SU(3)$-manifold!
%
\end{example}

\begin{example}\label{example;sun}
The latter observation generalizes to $K=\SU(n)$.  Let $T$ be the
standard maximal torus and $\alcove$ the alcove given in Example
\ref{example;centre-weyl} with $n$ vertices $\sigma_j$ corresponding
to the elements $c^j$ of the centre of $K$.  As invariant inner
product on $\k$ we take $(\xi,\eta)=-(4\pi^2)\inv\tr(\xi\eta)$.  This
choice of inner product makes the embedding $\t\to\R^n$ of Example
\ref{example;centre-weyl} an isometry, where we give $\t$ the induced
inner product and $\R^n$ the standard inner product.  Via the
isomorphism $\t\cong\t^*$ given by the inner product the vertex
$\sigma_j\in\t$ is identified with the $j$-th fundamental weight
$\varpi_j\in\t^*$.  The edge $\sigma_{01}$ of the alcove has
centralizer $K_{01}=\group{S}(\U(1)\times\U(n-1))$, so for $n>3$ the
corresponding stratum $X_{01}$ in $X$ consists of genuine
singularities.  Nevertheless we assert that its closure $\bar{X}_{01}$
possesses the structure of a smooth quasi-Hamiltonian
$\U(n)$-manifold.  Let $\sigma_{01}'=\{\,x\sigma_1\mid x>0\,\}$ be the
unique edge of the chamber $\ca{C}=\ca{C}\spcheck$ containing
$\sigma_{01}$.  Consider $\C^n$ with its standard Hermitian structure
$\langle z,w\rangle=z^*w$ and its standard $\U(n)$-action.  The
corresponding Hamiltonian $\U(n)$-structure $(\omega_0,\Phi_0)$ is
given by
$$
\omega_0=\Im\langle\cdot,\cdot\rangle=\frac{i}2\sum_jdz_j\,d\bar{z}_j,
\quad \Phi_0(z)=-2\pi^2i\,zz^*,
$$
where we have used the inner product to write $\Phi$ as a map into
$\k$.  Define $\ca{F}\colon K\times\sigma_{01}'\to\C^n$ by
$\ca{F}(k,x\sigma_1)=\sqrt{x/\pi}\,k_1$, where $k_1\in S^{2n-1}$
denotes the first column of the unitary matrix $k$.  The image of
$\ca{F}$ is $\C^n-\{0\}$.  Observe that $k_1=kv_1$, where
$v_1=(1,0,\dots,0)$ is the fundamental weight vector in $\C^n$, and
that the stabilizer of $v_1$ is $[K_{01},K_{01}]\cong\SU(n-1)$.  Hence
$\ca{F}$ descends to an injective map $f\colon X'_{01}\to\C^n$, where
$X'_{01}=K/[K_{01},K_{01}]\times\sigma_{01}'$.  Since $X'_{01}$ is
exactly the stratum of $X'=T^*K\impl$ associated with the face
$\sigma_{01}'$, it follows from \cite[Proposition
6.8(ii)]{guillemin-jeffrey-sjamaar} that $f$ is a symplectic
embedding.  (In fact $f$ is the restriction to a single stratum of a
symplectic embedding of $X'$ into the sum of the fundamental
representations.)  Clearly $f$ extends to a homeomorphism
$f\colon\bar{X}'_{01}=X'_{01}\cup X'_0\to\C^n$, which we regard as a
chart defining a smooth structure on $\bar{X}'_{01}$.  Pulling back
the symplectic structure on $\C^n$ via $f$ we obtain a symplectic
structure on $\bar{X}'_{01}$ for which $X'_{01}$ and $X'_0$ are smooth
symplectic submanifolds.  The chart $f$ is $K\times T$-equivariant
with respect to the (non-effective) action $(k,t)\cdot z=kt_1\inv z$
on $\C^n$, where $t=\diag(t_1,t_2,\dots,t_n)$, so $\bar{X}'_{01}$ is
isomorphic to $\C^n$ as a Hamiltonian $K\times T$-manifold.  There are
similar charts on $\bar{X}_{01}=X_{01}\cup X_0\cup X_1$.  We have
$\ca{F}(1,\sigma_1)=\pi^{-1/2}v_1$, so precomposing $f$ with $h\inv$,
where $h$ is the homeomorphism \eqref{equation;homeo}, we obtain a
homeomorphism
$$\phi_0\colon X_{01}\cup X_0\to D$$
where $D$ is the open ball of radius $1/\sqrt{\pi}$ in $\C^n$.  Using
$\phi_0$ as a chart we conclude from Proposition
\ref{proposition;cotangent-double} that $X_{01}\cup X_0$ is a smooth
quasi-Hamiltonian $K\times T$-manifold isomorphic to the
exponentiation of the Hamiltonian $K$-manifold $(D,\omega_0,\Phi_0)$.
The inverse $\phi_0\inv\colon D-\{0\}\to X_{01}$ is given by
\begin{equation}\label{equation;inverse}
z\mapsto\bigl(k[K_{01},K_{01}],\exp(\pi\norm{z}^2\sigma_1)\bigr),
\end{equation}
where $k\in K$ is any matrix whose first column is $z/\norm{z}$.
There is a similar chart
$$\phi_1\colon X_{01}\cup X_1\to D,$$
which is compatible with $\phi_0$.  It is constructed by precomposing
$\phi_0$ with a certain discrete symmetry of $X$ of the kind
considered in Propositions \ref{proposition;centre-action} and
\ref{proposition;duality}.  From Example \ref{example;centre-weyl} we
see that the action of the central element $c=e^{-2\pi i/n}I$ on
$\alcove$ sends $\sigma_j$ to $\sigma_{j+1}$, while the longest Weyl
group element $w_0$ sends $\sigma_j$ to $\sigma_{n-j}$.  (Here the
subscripts are taken modulo $n$.)  Therefore the action of $c\inv$
followed by that of $w_0$ is a symmetry of the alcove that maps the
edge $\sigma_{01}$ to itself, interchanging the two vertices
$\sigma_0$ and $\sigma_1$.  This Euclidean symmetry of $\alcove$ lifts
to the diffeomorphism $\eu{B}_{n_0}\circ\eu{A}_{\hat{c}}\inv\colon
DK\to DK$, where $\hat{c}$ is as in Example \ref{example;exact-centre}
and $n_0$ as in Example \ref{example;exact-long}.  The induced map on
the implosion $X=DK\impl$ maps $X_{01}$ to itself and interchanges
$X_0$ and $X_1$.  Since $\Ad\hat{c}(\sigma_1)=-\sigma_{n-1}$ we have
$\Ad\hat{c}\inv(\exp(x\sigma_1))=\exp(-x\sigma_{n-1})$ for $0<x<1$,
and so
\begin{multline*}
\pi(\hat{c})\inv\Ad\hat{c}\inv(\exp(x\sigma_1))
=c\inv\exp(-x\sigma_{n-1})\\
=\exp\sigma_{n-1}\exp(-x\sigma_{n-1})=\exp\bigl((1-x)\sigma_{n-1}\bigr).
\end{multline*}
Therefore $\eu{A}_{\hat{c}}\inv\bigl(k,\exp(x\sigma_1)\bigr)
=\bigl(k\hat{c},\exp((1-x)\sigma_{n-1})\bigr)$ and
$$
\eu{B}_{n_0}\circ\eu{A}_{\hat{c}}\inv\bigl(k,\exp(x\sigma_1)\bigr)
=\bigl(k\hat{c}n_0\inv,\exp((1-x)\sigma_1)\bigr).
$$
Hence the induced map $\psi\colon\bar{X}_{01}\to\bar{X}_{01}$ is given
by
$$
\psi\bigl((k,\exp(x\sigma_1)\mod{\sim}\bigr)
=\bigl(k\hat{c}n_0\inv,\exp((1-x)\sigma_1)\bigr)\mod{\sim}.
$$
Our second chart on $\bar{X}_{01}$ is then
$\phi_1=\phi_0\circ\psi\colon X_{01}\cup X_1\to D$.  An easy
calculation shows that the first column of the matrix
$k\hat{c}n_0\inv$ is $-k_1$.  Using \eqref{equation;inverse} we then
find that the transition map
$\phi_1\circ\phi_0\inv=\phi_0\circ\psi\circ\phi_0\inv$ is given by
\begin{equation}\label{equation;transition2n}
z\mapsto-\frac{\sqrt{\pi\inv-\norm{z}^2}}{\norm{z}}z.
\end{equation}
Therefore the pair of charts $(\phi_0,\phi_1)$ is a smooth atlas and
$\bar{X}_{01}$ is diffeomorphic to $S^{2n}$.  By Propositions
\ref{proposition;centre-action} and \ref{proposition;duality} the map
$\eu{B}_{n_0}\circ\eu{A}_{\hat{c}}\inv$ reverses the two-form on $DK$.
Consequently $\psi$ reverses the two-form on $\bar{X}_{01}$.  The
homomorphism $p\colon K\times T\to\U(n)$ defined by $p(k,t)=k(t_1\inv
I_n)$ is surjective; its kernel is the kernel of the $K\times
T$-action on $\bar{X}_{01}$; and the charts $\phi_0$ and $\phi_1$ are
$(K\times T)/\ker p=\U(n)$-equivariant.  Therefore $\bar{X}_{01}$ is
in effect a quasi-Hamiltonian $\U(n)$-manifold.  Comparing
\eqref{equation;transition2n} with the transition formula
\eqref{equation;transition} for the quasi-Hamiltonian $2n$-sphere
constructed in Appendix \ref{section;sphere} we draw the following
conclusion.

\begin{theorem}\label{theorem;s2n}
The closure of the stratum $X_{01}$ of $X=D\SU(n)\impl$ is a smooth
quasi-Hamiltonian $\U(n)$-manifold isomorphic to the spinning
$2n$-sphere of Theorem \ref{theorem;glue-disc}.  Under this
isomorphism the antipodal map on $S^{2n}$ corresponds to the
involution of $\bar{X}_{01}$ obtained by lifting the symmetry of the
alcove $\alcove$ that reverses the edge $\sigma_{01}$.
\end{theorem}

\end{example}

\begin{example}\label{example;u2h}
Let $K=\U(2,\H)=\Sp(2)$, the group of unitary $2\times2$-matrices over
the quaternions $\H$, with maximal torus $T=\bigl\{\diag(e^{2\pi
ix_1},e^{2\pi ix_2})\bigr\}\cong\U(1)^2$.  Identifying $\t$ with
$\R^2$ via the map $x\mapsto2\pi i\diag(x_1,x_2)$ we can write the
simple roots as
$$
(2\pi i)\inv\alpha_1(x)=x_1-x_2\qquad\text{and}\qquad(2\pi
i)\inv\alpha_2(x)=2x_2
$$
and the minimal root as $(2\pi i)\inv\alpha_0(x)=-2x_1$.  The alcove
is the isosceles right triangle $1/2>x_1>x_2>0$.  There are two
vertices, $\sigma_0$ resp.\ $\sigma_1$, corresponding to the central
elements $I$ resp.\ $-I$, and a third vertex $\sigma_2$ corresponding
to the element $\diag(-1,1)$.  The edge $\sigma_{01}$ is the
hypotenuse and the edges $\sigma_{02}$ and $\sigma_{12}$ are the right
edges.  Under the exponential map these correspond to elements of the
torus of the form $\diag(t,t)$, $\diag(t,1)$ and $\diag(-1,t)$,
respectively.
$$
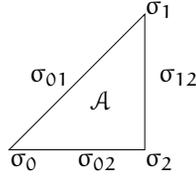
\setlength{\unitlength}{0.003mm}
\begin{picture}(7258,6813)(0,-10)
\path(12,654)(6012,654)(6012,6654)(12,654)(12,654)
\put(3612,2454){\makebox(0,0)[lb]{\smash{$\alcove$}}}
\put(12,-100){\makebox(0,0)[lb]{\smash{$\sigma_0$}}}
\put(6012,6654){\makebox(0,0)[lb]{\smash{$\sigma_1$}}}
\put(6012,-100){\makebox(0,0)[lb]{\smash{$\sigma_2$}}}
\put(3012,-100){\makebox(0,0)[lb]{\smash{$\sigma_{02}$}}}
\put(6612,3654){\makebox(0,0)[lb]{\smash{$\sigma_{12}$}}}
\put(912,3654){\makebox(0,0)[lb]{\smash{$\sigma_{01}$}}}
\end{picture}
$$
The centralizers and their commutator subgroups are as follows.
$$
\begin{tabular}{|l|l|l|}
\hline
$\sigma$&$K_\sigma$&$[K_\sigma,K_\sigma]$\\\hline\hline
$\alcove$&$T$&$1$\\
$01$&$\U(2)$&$\SU(2)$\\
$02$&$\U(1)\times\U(1,\H)$&$1\times\U(1,\H)$\\
$12$&$\U(1,\H)\times\U(1)$&$\U(1,\H)\times1$\\
$0$ or $1$&$K$&$K$\\
$2$&$\U(1,\H)\times\U(1,\H)$&$\U(1,\H)\times\U(1,\H)$\\\hline
\end{tabular}
$$
Here $\U(2)$ is the collection of matrices in $K$ with complex
coefficients, $\U(1)\times\U(1,\H)$ is the collection of matrices of
the form $\diag(a,b)$ with $a\in\C$, $b\in\H$ and $\abs{a}=\abs{b}=1$,
etc.  The alcove corresponds to a twelve-dimensional stratum
$\U(2,\H)\times\alcove$.  The strata corresponding to all three edges
are eight-dimensional and by Theorem \ref{theorem;smooth} are
removable singularities.  The strata $X_0$ and $X_1$ are genuine
singularities.  The stratum $X_2$ is diffeomorphic to
$K/K_2=\H\P^1\cong S^4$.  Note that $K_2$ is semisimple and is not the
centralizer of any element in the Lie algebra of $K$.  We do not know
of any relationship between the quasi-Hamiltonian structure of $S^4$
considered as an $\U(2,\H)$-conjugacy class and the
$\SU(2)\times\U(1)$-structure of Theorem \ref{theorem;s2n}.
As in Example \ref{example;sun} it can be shown that $\bar{X}_{02}$
and $\bar{X}_{12}$ are smooth quasi-Hamiltonian $\U(2,\H)$-manifolds
diffeomorphic to $\H\P^2$.
\end{example}

\section{Implosion and weighted projective spaces}
\label{section;projective}

\subsection*{The Hamiltonian case} 

The construction of an imploded cross-section as a subquotient of a
smooth manifold is satisfying in some respects, but less so in others.
In this section we give a different description of the universal
imploded cross-section $DK\impl$, namely as the orbit under $K$ of a
toric variety embedded in a representation of $K$.  In \cite[Section
6]{guillemin-jeffrey-sjamaar} a similar construction is performed in
the Hamiltonian case, which we review first.

Let $G=K^\C$\glossary{G@$G$, complexification of $K$} be the
complexification of $K$ and let $N$\glossary{N@$N$, maximal unipotent
subgroup of $G$} be the maximal unipotent subgroup of $G$ with Lie
algebra $\n=\bigoplus_{\alpha\in R_+}\g_\alpha$ (where $R_+$ denotes
the set of positive roots).  Let $\varpi_1$,\dots, $\varpi_r$ be the
fundamental weights of $G$, let $V_p$ denote the irreducible
$G$-module with highest weight $\varpi_p$, and define
$E=\bigoplus_{p=1}^rV_p$.  Let $E^N$ be the subspace of
$N$-invariants; this is
\begin{equation}\label{equation;en}
E^N=\bigoplus_{p=1}^r\C v_p,
\end{equation}
where $v_p\in V_p$ is a highest-weight vector.  The subset $G_N=GE^N$
of $E$ consists of finitely many $G$-orbits and is a closed affine
subvariety of $E$ (because it is equal to the closure of the $G$-orbit
through $v_1+v_2+\cdots+v_r$).  The module $E$ has a unique
$K$-invariant Hermitian inner product such that each $v_p$ has length
$1$; and we equip the orbits in $G_N$ with the K\"ahler structures
induced by this flat K\"ahler structure on $E$.  In \cite[Section
6]{guillemin-jeffrey-sjamaar} is constructed a $K\times T$-equivariant
homeomorphism $T^*K\impl\to G_N$, which maps each stratum in
$T^*K\impl$ symplectomorphically onto a $G$-orbit in $G_N$.  This
allows us to view $T^*K\impl$ as a complex algebraic variety with a
K\"ahler structure on each stratum.  Observe also that $G_N=KE^N$
because of the Iwasawa decomposition $G=KAN$ (with $A=\exp i\t$).  In
addition, $E^N$ with its natural $T$-action and moment map
$\Phi\bigl(\sum_pz_pv_p\bigr)=\sum_p\abs{z_p}^2\varpi_p$ is the
symplectic toric variety associated with the polyhedron
$\bar{\ca{C}}$.  (See \cite[Remark 6.7]{guillemin-jeffrey-sjamaar},
but note the different normalization of the moment map.)

It is clear from the examples in Section \ref{section;double} that the
quasi-Hamiltonian imploded cross-section $DK\impl$ is \emph{not} a
complex algebraic variety.  (Indeed $DK\impl=S^4$ for $K=\SU(2)$.)
Nevertheless it turns out that we can still build $DK\impl$ as the
$K$-orbit of a toric variety embedded in a $K$-module.

\subsection*{Construction of a toric variety}

To adapt the construction of $G_N$ to the quasi-Ham\-ilton\-ian
category, we start by replacing the chamber $\ca{C}$ by the alcove
$\alcove$ and building a toric variety $X$ corresponding to it.  (We
use the isomorphism $\t\to\t^*$ given by the inner product to view
$\alcove$ as a subset of $\t^*$.)  We may assume without loss of
generality that $K$ is almost simple.  Then the alcove is a simplex,
so the toric variety is just a weighted projective space and can be
characterized as a symplectic cut of the $T$-module $E^N$ considered
in \eqref{equation;en}.  Let $\alpha_1$,\dots, $\alpha_r$ be the
simple roots of $K$ and $\alpha_{r+1}$ the highest root.  Let $H$ be
the circle subgroup of $T$ generated by $2\pi i\,\alpha_{r+1}^\vee$
and let $X$ be the symplectic cut (in the sense of
\cite{lerman;symplectic-cuts;mathematical-research}) of $E^N$ relative
to $H$ at level $k=2/\inner{\alpha_{r+1},\alpha_{r+1}}$.  (Modulo the
identification $\t\cong\t^*$ we have
$\alpha_{r+1}^\vee=2\alpha_{r+1}/\inner{\alpha_{r+1},\alpha_{r+1}}$,
so $\alpha_{r+1}^\vee=k\alpha_{r+1}$.)  Then $X$ is a symplectic
orbifold which inherits a Hamiltonian $T$-action from $E^N$ with a
moment map which we shall also denote by $\Phi$.  The moment polytope
of $X$ is obtained by truncating the moment polytope of $E^N$ at the
hyperplane $\lambda(\alpha_{r+1}^\vee)=k$, i.e.\
\begin{equation}\label{equation;cut}
\begin{split}
\Phi(X)&=\Phi(E^N)\cap\bigl\{\,\lambda\in\t^*\bigm|(2\pi
i)\inv\lambda(\alpha_{r+1}^\vee)\le k\,\bigr\}\\
&=\bar{\ca{C}}\cap\bigl\{\,\lambda\in\t^*\bigm|(2\pi
i)\inv\inner{\lambda,\alpha_{r+1}}\le1\,\bigr\}\\
&=\bar{\alcove}.
\end{split}
\end{equation}
Note that the size of the image of $\alcove$ in $\t^*$, and hence the
symplectic structure on $X$, depend on the choice of the inner product
on $\k$.  The $T$-action on $E^N$ is effective and therefore so is the
induced $T$-action on $X$.  We conclude that $X$ is a symplectic toric
orbifold in the sense of
\cite{lerman-tolman;hamiltonian-torus-actions-symplectic-orbifolds}.
By [\emph{loc.\ cit.}, Theorems 7.4, 8.1] such orbifolds are
classified by \emph{labelled polytopes}, simple rational polytopes in
$\t^*$ with a label attached to each face of codimension one.  The
label is the integer $l$ such that every element in the preimage of
the face under the moment map has orbifold structure group $\Z/l\Z$.
We can now describe $X$ more explicitly as follows.  The coefficients
$m_1$,\dots, $m_r$ in the expansion of the coroot of the highest root,
$$
\alpha_{r+1}^\vee=\sum_{p=1}^rm_p\alpha_p^\vee,
$$
are relatively prime positive integers.  Set $m_{r+1}=1$ and
$m=(m_1,\dots,m_{r+1})$.

\begin{lemma}
\begin{enumerate}
\item\label{item;weight}
$X\cong\C\P^r_m$, the weighted projective $r$-space with weights $m$.
\item\label{item;label}
$X$ is the symplectic toric orbifold associated with the polytope
$\bar\alcove$ with label $1$ attached to every codimension-one face.
\item\label{item;stabilizer}
Let $x\in X$ and $\sigma\le\alcove$.  Suppose $\Phi(x)\in\sigma$.
Then $T_x=T\cap[K_\sigma,K_\sigma]$.
\end{enumerate}
\end{lemma}

\begin{proof}
The symplectic cut $X$ is by definition the symplectic quotient of
$E^N\times\C$ with respect to $H$ at level $k$.  Here the $H$-action
on $E^N$ is the restriction of the $T$-action and the $H$-action on
$\C$ is defined infinitesimally by $\alpha_{r+1}^\vee\cdot w=w$.
Introducing coordinates $(z_1,\dots,z_r,z_{r+1})$ on $E\times\C$
(relative to the basis $v_1$,\dots, $v_r$ of $E$) and writing
$h=\exp(t\alpha^\vee_{r+1})$ we have
$$
h\cdot v=\bigl(e^{2\pi it\varpi_1(\alpha^\vee_{r+1})}z_1,\dots,
e^{2\pi it\varpi_r(\alpha^\vee_{r+1})}z_r\bigr)=\bigl(e^{2\pi
im_1t}z_1,\dots, e^{2\pi im_rt}z_r\bigr),
$$
since $\varpi_q(\alpha^\vee_p)=\delta_{pq}$.  Thus $X$ is the
symplectic quotient of $\C^{r+1}$ by the $\U(1)$-action with weights
$m$, which proves \eqref{item;weight}.  Moreover, by
\eqref{equation;cut} the $T$-moment map image of $X$ is the rational
simplex $\bar\alcove$.  To calculate the label attached to a face
$\sigma_p$, note that for $p\le r$ the $T$-moment map $\Phi$ maps (a
suitable multiple of) $x_p=\bigl(\sum_{q\ne p}v_q,1\bigr)$ into
$\sigma_p$.  The $H$-stabilizer of $x_p$ is trivial because $H$ acts
with weight $1$ on $\C$.  Thus $\sigma_p$ has label $1$ for $p\le r$.
Moreover, (a suitable multiple of) the point
$\bigl(\sum_{p=1}^rv_p,0\bigr)$ is mapped into $\sigma_{r+1}$ and this
point also has trivial $H$-stabilizer because $\gcd(m_1,\dots,m_r)=1$.
Thus $\sigma_{r+1}$ has label $1$.  Finally, by \cite[Lemma
6.6]{lerman-tolman;hamiltonian-torus-actions-symplectic-orbifolds} the
stabilizer $T_x$ for $x\in\Phi\inv(\sigma)$ is the torus whose Lie
algebra is the subspace of $\t$ spanned by all $\alpha^\vee$ that are
perpendicular to $\sigma$.  This shows that
$\t_x=\t\cap[\k_\sigma,\k_\sigma]$.  Now $T\cap[K_\sigma,K_\sigma]$ is
a maximal torus of $[K_\sigma,K_\sigma]$.  In particular it is
connected and therefore $T_x=T\cap[K_\sigma,K_\sigma]$.
\end{proof}

\begin{example}\label{example;planches}
Using \cite[Ch.~6, Planches I-IX]{bourbaki;groupes-algebres} one
easily calculates the following table.
$$
\begin{tabular}{|l|c|}
\hline
$K$&$m$\\\hline\hline
$\A_r$&$(1,1,\dots,1)$\\
$\B_r$&$(1,2,2,\dots,2,1,1)$\\
$\C_r$&$(1,1,\dots,1)$\\
$\D_r$&$(1,2,2,\dots,2,1,1,1)$\\
$\E_6$&$(1,2,2,3,2,1,1)$\\
$\E_7$&$(2,2,3,4,3,2,1,1)$\\
$\E_8$&$(2,3,4,6,5,4,3,2,1)$\\
$\F_4$&$(2,3,2,1,1)$\\
$\G_2$&$(1,2,1)$\\
\hline
\end{tabular}
$$
The upshot is that $X$ is the standard projective space $\C\P^r$ for
$K=\SU(r+1)$ or $K=\U(r,\H)$ and a weighted projective space in all
other cases.
\end{example}

\subsection*{Embedding}

We want to embed $X$ into a $K$-module $V$ in such a way that the
points in $\Phi\inv(\sigma)$ have stabilizer $[K_\sigma, K_\sigma]$;
this will then give us an identification $DK\impl\cong KX$.  For each
$\sigma\le\alcove$ one should choose a representation $W_\sigma$ and a
point $p_\sigma\in W_\sigma$ whose stabilizer is
$[K_\sigma,K_\sigma]$.

For the vertices, one can use the fact that
$K_\sigma=[K_\sigma,K_\sigma]$: one then just needs to take a faithful
representation of the group, which embeds the group into a vector
space $V$, and have $K$ act by the adjoint representation.

For the edges, we want, for each pair of roots $\alpha_i$, $\alpha_j$
generating the faces of the alcove, a representation $V_{i,j}$ and a
point $p_{i,j}$ such that the stabilizer $[K_{i,j},K_{i,j}]$ of
$p_{i,j}$ is generated by the root spaces $e_{\alpha_k}$,
$f_{\alpha_k}$, for $k\ne i,j$, and such that the torus $T$ acts on
the line generated by $p_{i,j}$ by the weight
$l_j\varpi_j-l_i\varpi_i$.  Here we define $l_i=m/m_i$, with
$m=\lcm(m_1,\dots,m_{r+1})$.

A first case is when either $i$ or $j$, say $i$, is $r+1$.  One then
just takes the $l_j$-th power of the $j$-th fundamental
representation, and $p_{i,j}$ the highest weight vector.

Suppose now that $i$ and $j$ both correspond to simple roots.  Then
\begin{equation}\label{equation;muvarpi}
\biginner{l_j\varpi_j-l_i\varpi_i,\alpha^\vee_k}
\begin{cases}
>0&\text{for $k=j$}\\
<0&\text{for $k=i$}\\
=0&\text{otherwise.}
\end{cases}
\end{equation}
(For the lowest root, this uses the definition of the $l_i$.)  We can
use the Weyl group to bring $l_j\varpi_j-l_i\varpi_i$ into the chamber
$\ca{C}$, and the roots $\alpha_k$ for $k\ne i$ into the set of
positive roots.  Because of \eqref{equation;muvarpi}, the
representation with highest weight $l_j\varpi_j-l_i\varpi_i$ will then
be a multiple of one of the highest weight representations.

Now to the embedding.  We proceed in an inductive fashion on the faces
of the alcove, starting with the vertices.  Let $(a_1,\dots,a_{r+1})$
denote the projective coordinates of a point in the weighted
projective space $X$, where we can assume that $(a_1,\dots,a_{r+1})$
is of norm one.

The vertices correspond to points where only one $a_i$ is non-zero,
and we map this to $W_{i,i}$, the representation corresponding to the
vertex, by
$$a_i\mapsto(\bar a_i)^{l_i}(a_i)^{l_i}p_{i,i}$$
Now for the edges.  For each pair $(i,j)$ with $1\le i,j\le r+1$
consider the $(i,j)$-th edge $\sigma=\sigma_{i,j}$, and choose a
representation $W_{i,j}$ and a point $p_{i,j}$ with stabilizer
$[K_{\sigma_{i,j}},K_ {\sigma_{i,j}}]$.  (Clearly we may assume that
$W_{i,j}=W_{j,i}$.)  We map the edge by
$$(a_i, a_j)\mapsto(\bar a_i)^{l_i}(a_j)^{l_j}p_{i,j}.$$
This is sufficient by Proposition
\ref{proposition;centralizer}\eqref{item;centralizer-commutator}.  In
other words, we extend this map to a map from $X$ to
$V=\bigoplus_{i,j=1}^{r+1}W_{i,j}$ by
$$
(a_1,\dots,a_{r+1})\mapsto\bigoplus_{i,j=1}^{r+1}(\bar a_i)^{l_i}
(a_j)^{l_j} p_{i,j}.
$$
(Note that for $i\ne j$ the summand $W_{i,j}=W_{j,i}$ occurs twice in
$V$.)  We then have the desired embedding of $X$, and a homeomorphism
$DK\impl\cong KX$.

\begin{example}
For $\SU(r+1)$ we can proceed slightly differently and embed
$X=\C\P^r$ into a sum of weight spaces as follows.  Consider the unit
sphere
$$\bigl\{\,(a_1,\dots,a_{r+1})\bigm|\abs{a_1}^2+\cdots+\abs{a_{r+1}}^2
=1\,\bigr\}$$
in $\C^{r+1}$.  As the vertices of the alcove all correspond to
central elements, we can set $W_{i,i}=\C$, $p_{i,i}=1$.  For the
edges, if $e_1$,\dots, $e_{r+1}$ denotes the standard basis of
$\C^{r+1}$, we set $ p_{i,j}=e_{i+1}\wedge\cdots\wedge e_j$ for $i<j$,
and $E_{i,j}=e_{i+1}\wedge\cdots\wedge e_{r+1}\wedge
e_1\wedge\cdots\wedge e_j$ for $i>j$, where these elements live in the
representations $W_{i,j}=\Lambda^{j-i}\C^{r+1}$ for $i\le j$, and
$W_{i,j}=\Lambda^{r+1+j-i}\C^{r+1}$ for $i>j$.  Set $V_{i,j}=\C
E_{i,j}$.  Map the unit sphere to $U=\bigoplus V_{i,j}$ by
$$(c_1,\dots,c_{r+1})\mapsto\sum_{i,j=1}^{r+1}\bar{c}_ic_jE_{i,j}.$$
This map is $T$-equivariant, and factors through $X$.  We set
$W=\bigoplus_{i,j=1}^{r+1}W_{i,j}$.  The orbit $\SU(r+1)\cdot X$ is
then our universal imploded cross-section.  If $(a_1,\dots,a_{r+1})$
are the coordinates of a point $p$ in $X$, with $a_{i_1}$,
$a_{i_2}$,\dots, $a_{i_s}$ the non-zero coordinates (with
$i_1<i_2<\cdots<i_s$), the stabilizer of $p$ is
$\SU(r+1+i_1-i_s)\times\SU(i_2-i_1)\times\cdots\times\SU(i_s-i_{s-1})$.
This is the commutator subgroup of the stabilizer under the adjoint
action of elements
$$
\diag(l_1,\dots,l_1,l_2,\dots,l_2,l_3,\dots,l_3,\dots,
l_{s-1},l_s,\dots,l_s),
$$
where $l_1$ is repeated $r+1+i_1-i_{s}$ times, and the other $l_j$
are repeated $i_j-i_{j-1}$ times.
\end{example}

\section{The master moduli space}\label{section;parabolic}

Let $\Sigma$ be a compact oriented two-manifold of genus $g$ with
boundary components $V_1$, $V_2$,\dots, $V_n$.  As before, let $K$ be
a compact simply connected Lie group and choose an $n$-tuple
$C=(C_1,C_2,\dots,C_n)$ of conjugacy classes in $K$.  Consider
homomorphisms $\pi_1(\Sigma)\to K$ which map the cycle $[V_i]$ into
the class $C_i$ for $1\le i\le n$.  The space $M(\Sigma,C)$ of all
such homomorphisms modulo conjugation by $K$ can be interpreted as the
space of flat connections on the trivial principal bundle
$\Sigma\times K$ with prescribed holonomy $c_i$ around $V_i$, modulo
gauge equivalence.  Atiyah and Bott
\cite{atiyah-bott;yang-mills-equations-riemann} proved that
$M(\Sigma,C)$ carries a natural symplectic structure on a dense open
subset.  Alekseev et al.\ \cite[Theorem
9.2]{alekseev-malkin-meinrenken;lie-group} proved that, for $n\ge1$,
$M(\Sigma,C)$ is a quasi-symplectic quotient,
\begin{equation}\label{equation;moduli-quasi}
M(\Sigma,C)\cong M(\Sigma)\quot[c]K^n.
\end{equation}
(On the right we reduce at the unique point $c$ where $C$ intersects
$\bigl(\exp\bar{\alcove}\bigr)^n$.)  Here $M(\Sigma)$ is the
quasi-Hamiltonian $K^n$-manifold
\begin{equation}\label{equation;moduli-fusion}
M(\Sigma)=\D K^{\fusion g}\fusion DK^{\fusion(n-1)};
\end{equation}
$\D K$ denotes the quasi-Hamiltonian $K$-manifold obtained from $DK$
by fusing the two $K$-actions, and $X^{\fusion l}$ denotes the
$l$-fold fusion product of a space $X$ with itself.  Thus
$M(\Sigma)=K^{2(g+n-1)}$ as a manifold.  The isomorphism
\eqref{equation;moduli-quasi} implies, by Theorem
\ref{theorem;singular-quotient}, that $M(\Sigma,C)$ is a stratified
space with symplectic strata, thus confirming the Atiyah-Bott theorem.
The definition of the isomorphism involves a choice of base points
$p_i\in V_i$ for $1\le i\le n$, which enables one to identify
$M(\Sigma)$ with the space of flat connections on $\Sigma\times K$
modulo gauge transformations that are the identity at the $p_i$.  (See
\cite[Theorems 9.1, 9.3]{alekseev-malkin-meinrenken;lie-group}.)

The first goal of this section is to build a \emph{master moduli
space} $M$, which is symplectic and from which all moduli spaces
$M(\Sigma,C)$ can be obtained by symplectic reduction with respect to
a compact torus.  Given the results of Section \ref{section;implosion}
this is entirely straightforward: we simply define $M=M(\Sigma)\impl$,
the implosion of $M(\Sigma)$ with respect to $K^n$.  This carries an
action of the $rn$-dimensional torus $T^n$ (where $r$ is the rank of
$K$) and an imploded moment map $\Phi\colon M\to T^n$ taking values in
$\bigl(\exp\bar{\alcove}\bigr)^n$.  Theorem
\ref{theorem;decomposition} and Addendum \ref{addendum;quotient}
immediately imply the following result.

\begin{theorem}\label{theorem;master}
The master moduli space $M=M(\Sigma)\impl$ has a $T^n$-invariant
locally finite decomposition into symplectic manifolds.  The
$T^n$-action on each stratum is quasi-Hamiltonian with moment map
given by the restriction of $\Phi$.  Let $C$ be a conjugacy class in
$K^n$ and $c$ its unique intersection point with
$\bigl(\exp\bar{\alcove}\bigr)^n$.  Then $M\quot[c]T^n\cong
M(\Sigma,C)$.
\end{theorem}

Let $\Sigma'$ be the closed Riemann surface of genus $g$ obtained by
capping off the boundary components of $\Sigma$ and choose a complex
structure on $\Sigma'$.  A well-known theorem due to Mehta and
Seshadri \cite{mehta-seshadri} establishes an equivalence between the
moduli space $M(\Sigma,C)$ and the moduli space of holomorphic vector
bundles on $\Sigma'$ with prescribed parabolic structures at the
points $p_1$, $p_2$,\dots, $p_n$.  In
\cite{hurtubise-jeffrey-sjamaar;moduli-framed} we will prove, for
$K=\SU(n)$, an analogue of the Mehta-Seshadri theorem by showing that
the master moduli space $M$ possesses a complex algebraic structure
such that every moduli space of parabolic bundles can be written as a
GIT quotient of $M$ with respect to the complex torus $(T^\C)^n$ (and
we will indicate a possible approach for arbitary $K$).

We now want to give a more explicit description of $M$ and also draw a
parallel between the various moduli spaces associated to $\Sigma$ on
one hand and the following diagram of spaces on the other hand:
\begin{equation}\label{equation;doubles}
\vcenter{\xymatrix{DK\fusion DK\ar@{.>}[r]^-{\quot[1]K}
&DK\ar@{.>}[dr]_{\quot[c_0]K}&
K\times\exp\bar{\alcove}\ar@{))->}[l]\ar@{>>}[r]&
DK\impl\ar[r]^-{(\Psi_1)\impl}\ar@{.>}[dl]^{\quot[c_0]T}&K\\
&&C_0}}
\end{equation}
Here the dotted arrows indicate quasi-Hamiltonian reduction with
respect to the relevant group, $C_0$ denotes a conjugacy class of $K$
and $c_0$ its unique intersection point with $\exp\bar{\alcove}$,
$K\times\exp\bar{\alcove}$ is the inverse image of $\exp\bar{\alcove}$
in $DK$, and $K\times\exp\bar{\alcove}\to DK\impl$ is the quotient map
for the equivalence relation $\sim$.  The moduli spaces analogous to
$DK$, $DK\impl$ and $C_0$ are $M(\Sigma)$, $M=M(\Sigma)\impl$ and
$M(\Sigma,C)$, respectively.  To find the missing spaces, we punch an
extra hole with boundary circle $V_{n+1}$ in $\Sigma$ to obtain a new
Riemann surface $\hat{\Sigma}$.  Because of
\eqref{equation;moduli-fusion} the moduli space of $\hat{\Sigma}$ is
then the quasi-Hamiltonian $K^{n+1}$-manifold
$$
M(\hat{\Sigma})=\D K^{\fusion g}\fusion DK^{\fusion n},
$$
and $M(\Sigma)=M(\hat{\Sigma})\quot[1]K$ by Lemma
\ref{lemma;universal}\eqref{item;universal-centre}.  (Indeed, if for
$n=0$ we \emph{define} $M(\Sigma)=M(\hat{\Sigma})\quot[1]K$, then
\eqref{equation;moduli-quasi} is true even for $n=0$.)  Thus
$M(\hat{\Sigma})$ is the moduli space analogue of the space $DK\fusion
DK$ on the left in \eqref{equation;doubles}.  According to
\cite[\S~9.2]{alekseev-malkin-meinrenken;lie-group}, in suitable
coordinates $(a,b,u,v)\in K^g\times K^g\times K^n\times K^n$ on
$M(\hat{\Sigma})$ the action of $(k_1,k_2,\dots,k_{n+1})\in K^{n+1}$
is given by
\begin{alignat*}{2}
a_h&\mapsto\Ad(k_{n+1})a_h,&\qquad b_h&\mapsto\Ad(k_{n+1})b_h,\\
u_m&\mapsto k_{n+1}u_mk_m\inv,&\qquad v_m&\mapsto\Ad(k_m)v_m,
\end{alignat*}
and the moment map $\Phi=(\Phi_1,\Phi_2,\dots,\Phi_{n+1})$ by
\begin{equation}\label{equation;commutator}
\begin{split}
\Phi_i(a,b,u,v)&=v_i\qquad\text{for $i=1$, $2$,\dots, $n$,}\\
\Phi_{n+1}(a,b,u,v)
&=\prod_{h=1}^g[a_h,b_h]\prod_{m=1}^n\Ad(u_m)v_m\inv.
\end{split}
\end{equation}
The $(a,b,u,v)$ are to be thought of as holonomies along paths $A_h$,
$B_h$ ($1\le h\le g$) and $U_m$, $V_m$ ($1\le m\le n$).  The $A_h$ and
$B_h$ originate at a base point $p_{n+1}$ on the boundary component
$V_{n+1}$ of $\hat{\Sigma}$ and represent the standard homology basis
of the closed surface $\Sigma'$, and the $U_m$ are paths joining
$p_{n+1}$ to $p_m$, the base points on the boundary components $V_m$
of $\Sigma$.  We can further assume that these paths do not intersect
except at their endpoints, and are such that we can write $\Sigma$ in
the usual way as the quotient of a polygon with sides (in order)
\begin{multline*}
A_1,B_1,A_1\inv,B_1\inv,A_2,B_2,A_2\inv,B_2\inv,\dots,
A_g,B_g,A_g\inv,B_g\inv,\\
U_1,V_1\inv,U_1\inv,U_2,V_2\inv,U_2\inv,\dots,U_n,V_n\inv,U_n\inv,
\end{multline*}
with identifications as given along the edges.  The fundamental group
of $\Sigma$ is then the group generated by $A_h$, $B_h$, $U_mV_m\inv
U_m\inv$ subject to the relation
$$\prod_{h=1}^g[A_h,B_h]\prod_{m=1}^nU_mV_m\inv U_m\inv=1.$$
Using the coordinate system $(a,b,u,v)$ we can describe $M(\Sigma)$ as
\begin{equation}\label{equation;moduli}
M(\Sigma)=M(\hat{\Sigma})\quot[1]K=\bigl\{\,(a,b,u,v)\in
K^{2(g+n)}\bigm|\Phi_{n+1}(a,b,u,v)=1\,\bigr\}\big/K,
\end{equation}
where $K$ acts as the last factor of $K^{n+1}$.  The preimage of
$\bigl(\exp\bar{\alcove}\bigr)^n$ under the moment map for the
residual $K^n$-action is then
$$
M(\Sigma,\alcove)=\bigl\{\,(a,b,u,v)\in
K^{2(g+n)}\bigm|\Phi_{n+1}(a,b,u,v)=1,\quad
v\in\bigl(\exp\bar{\alcove}\bigr)^n\,\bigr\}\big/K.
$$
This is to be viewed as the analogue of the space
$K\times\exp\bar\alcove$ in \eqref{equation;doubles} and can be
interpreted as the moduli space of flat connections whose holonomy
around the loops $V_i$ is restricted to lie in $\exp\bar\alcove$.
Implosion gives
$$
M=M(\Sigma)\impl=\bigl\{\,(a,b,u,v)\in
K^{2g}\times(DK\impl)^n\bigm|\Phi_{n+1}(a,b,u,v)=1\,\bigr\}\big/K,
$$
where $(u_m,v_m)\in K\times\exp\bar\alcove$ represents an element of
$DK\impl$.  From \eqref{equation;commutator} we see that the moment
map $M\to T^n$ is the map induced by the holonomy map
$(a,b,u,v)\mapsto v$.  Reduction at
$c\in\bigl(\exp\bar\alcove\bigr)^n$ gives
$$
M(\Sigma,C)=\bigl\{\,(a,b,u,c)\in K^{2g}\times(DK\impl/T)^n
\bigm|\Phi_{n+1}(a,b,u,c)=1\,\bigr\}\big/K,
$$
where $(u_m,c_m)\in K\times\exp\bar\alcove$ represents an element of
$DK\impl/T$.  Finally, setting $u=1$ in \eqref{equation;moduli} gives
the space of all representations of the fundamental group,
$$
\Hom(\pi_1(\Sigma),K)/K=\bigl\{\,(a,b,v)\in
K^{2g+n}\bigm|\Phi_{n+1}(a,b,1,v)=1\,\bigr\}\big/K,
$$
and a natural map $f\colon M(\Sigma)\to\Hom(\pi_1(\Sigma),K)/K$,
which descends to a map $f\impl$ from $M(\Sigma)\impl$ to
$\Hom(\pi_1(\Sigma),K)/K$.  The moduli space counterpart of diagram
\eqref{equation;doubles} is thus
$$
\xymatrix{M(\hat{\Sigma})\ar@{.>}[r]^-{\quot[1]K}
&M(\Sigma)\ar@{.>}[dr]_{\quot[c]K^n}&
M(\Sigma,\alcove)\ar@{))->}[l]\ar@{>>}[r]&
M(\Sigma)\impl\ar[r]^-{f\impl}\ar@{.>}[dl]^{\quot[c]T^n}
&\Hom(\pi_1(\Sigma),K)/K\\
&&M(\Sigma,C)}
$$
We do not know if $f$ or $f\impl$ can be interpreted as a moment map
in any sense.

\appendix

\section{The conjugation action}\label{section;adjoint}
 
This appendix is a review of some facts concerning the conjugation, or
adjoint, action of a compact Lie group on itself.  Several of these
facts are well-known, but for lack of a convenient reference we
provide short proofs.  Let $K$ be a simply connected compact Lie group
with Lie algebra $\k$.  Fix a maximal torus $T$ of $K$ and a chamber
$\ca{C}$ in the dual of the Cartan subalgebra $\t$.  Let
$\ca{C}\spcheck$ be the dual chamber in $\t$ and let $\alcove$ be the
unique alcove in $\ca{C}\spcheck$ such that the origin is in the
closure of $\alcove$.

\subsection*{Centralizers}

Let $K_g=\{\,h\in K\mid hg=gh\,\}$ denote the centralizer of $g\in K$.
 
\begin{lemma}\label{lemma;borel}
The centralizer $K_g$ is connected for all $g\in K$.  Let
$\sigma$\glossary{sigma@$\sigma$, face of chamber or alcove} be a
(relatively open) face of the alcove $\alcove$.  Then
$K_{\exp\xi}=K_{\exp\eta}$ for all $\xi$, $\eta\in\sigma$.
\end{lemma}

\begin{proof}
The first assertion is
\cite[Corollaire~3.5]{borel;sous-groupes-commutatifs}.  (Cf.\ also
\cite[Ch.~9, \S~5.3, Corollaire 1]{bourbaki;groupes-algebres}.)  Let
$\xi\in\sigma$ and put $g=\exp\xi$.  To prove the second assertion we
need only show that the Lie algebra $\k_g=\ker(\Ad(g)-1)$ of $K_g$ is
independent of the point $\xi\in\sigma$.  Let $R$\glossary{R@$R$, root
system of $(K,T)$} denote the root system of the pair $(K,T)$ and
$\g=\k^\C$ the complexification of $\k$.  Since $g$ is in the torus
$T$, $K_g$ contains $T$, so the root-space decomposition
$\g=\t^\C\oplus\bigoplus_{\alpha\in R}\g_\alpha$ restricts to the
root-space decomposition
$$
(\k_g)^\C=\t^\C\oplus\bigoplus_{\alpha\in R_\xi}\g_\alpha,
$$
where $R_\xi\subset R$ is the root system of $(K_g,T)$.  To determine
$R_\xi$ choose a root vector $x_\alpha\in\g_\alpha$ for each
$\alpha\in R$.  Then $x_\alpha\in(\k_g)^\C$ if and only if
$e^{\ad\xi}x_\alpha=\Ad(g)x_\alpha=x_\alpha$.  Identifying as usual a
root $\alpha\colon T\to\U(1)$ with its differential
$T_1\alpha\colon\t\to i\R$ we have
$[\xi,x_\alpha]=\alpha(\xi)x_\alpha$, so
$$
R_\xi=\{\,\alpha\in R\mid\alpha(\xi)\in2\pi i\Z\,\}.
$$
Now let
\begin{equation}\label{equation;root}
R_\sigma=\{\,\alpha\in R\mid\alpha(\eta)\in2\pi i\Z\text{ for all
$\eta\in\sigma$}\,\}.\glossary{R@$R_\sigma$, root system of
$(K_\sigma,T)$}
\end{equation}
To finish the proof it suffices to show that $R_\xi=R_\sigma$.
Clearly $R_\sigma\subset R_\xi$.  Let $\alpha\in R_\xi$.  Then
$\alpha(\xi)=2\pi in$ with $n\in\Z$, so $\xi$ is in the singular
hyperplane
$$
H_{\alpha,n}=\{\,\eta\in\t\mid\alpha(\eta)=2\pi
in\,\}.\glossary{H@$H_{\alpha,n}$, singular hyperplane of $\t$}
$$
The intersection of $H_{\alpha,n}$ with the closed alcove
$\bar{\alcove}$ is the closure of a face $\tau$ of $\alcove$.  Since
$\xi\in\sigma$, $\sigma$ intersects the closure of $\tau$ and
therefore is a face of $\tau$.  Hence $\sigma\subset H_{\alpha,n}$.
This implies $\alpha(\eta)=2\pi in$ for all $\eta\in\sigma$, so
$\alpha\in R_\sigma$.
\end{proof}

This result is well-known to be false if $K$ is not simply connected.
One need only think of a rotation by $\theta$ in $\SO(3)$.  If
$\theta$ is not $\pi$, the centralizer is $\SO(2)$; if it is $\pi$,
the centralizer is $\group{O}(2)$.

If $\sigma$ is the face containing $\xi\in\bar{\alcove}$ we shall
permit ourselves to call $K_{\exp\xi}$ the \emph{centralizer of the
face} $\sigma$.  We denote it by $K_\sigma$\glossary{K@$K_\sigma$,
centralizer of a face} and its Lie algebra by $\k_\sigma$.  Its root
system $R_\sigma$ is given by \eqref{equation;root}, so its
complexified Lie algebra is
$$
(\k_\sigma)^\C=\t^\C\oplus\bigoplus_{\alpha\in
R_\sigma}\g_\alpha.\glossary{K@$\k_\sigma$, Lie algebra of $K_\sigma$}
$$
The proof of Lemma \ref{lemma;borel} also yields the following useful
information.  Here $Z(G)$\glossary{Z@$Z(A)$, centre of a group $A$}
denotes the centre and $G^0$\glossary{A0@$A^0$, unit component of a
Lie group $A$} the identity component of a Lie group $G$, and $\z(\g)$
denotes the centre of a Lie algebra $\g$.  For simplicity we denote
the face $\{0\}$ of $\alcove$ by $0$.

\begin{lemma}\label{lemma;face}
Let $\sigma$ and $\tau$ be faces of the alcove $\alcove$.
\begin{enumerate}
\item\label{item;face}
If $\sigma\le\tau$, then $K_\tau\subset K_\sigma$.
\item\label{item;face-face}
If $\sigma\ge0$ and $\sigma'$ is the unique face of the chamber
$\ca{C}\spcheck$ containing $\sigma$, then $K_\sigma=K_{\sigma'}$, the
centralizer of $\sigma'$ for the adjoint action of $K$ on $\k$.
\item\label{item;face-centre}
Let $H_\sigma$ be the affine subspace of $\t$ spanned by $\sigma$.
Then the linear subspace of $\t$ parallel to $H_\sigma$ is
$\z(\k_\sigma)$, and $\exp H_\sigma$ is contained in $Z(K_\sigma)$.
If $\sigma\ge0$, then $H_\sigma=\z(\k_\sigma)$ and $\exp
H_\sigma=Z(K_\sigma)^0$.
\item\label{item;face-shift}
There exists $g_\sigma\in[K_\sigma,K_\sigma]\cap Z(K_\sigma)$, unique
up to multiplication by an element of
$\Gamma_\sigma=[K_\sigma,K_\sigma]\cap
Z(K_\sigma)^0$\glossary{Gamma@$\Gamma_\sigma$, covering group
$[K_\sigma,K_\sigma]\cap Z(K_\sigma)^0$}, such that $\exp
H_\sigma\subset g_\sigma Z(K_\sigma)^0$.
\end{enumerate}
\end{lemma}

\begin{proof}
By \eqref{equation;root} $\sigma\le\tau$ implies $R_\tau\subset
R_\sigma$, and hence $K_\tau\subset K_\sigma$.

Assume $\sigma\ge0$ and let $\alpha\in R_\sigma$.  Then
$\alpha(\xi)=0$ for any $\xi\in\sigma$.  Letting $\alpha_1$,
$\alpha_2$,\dots, $\alpha_r$ be the simple roots with respect to the
chamber $\ca{C}$ and writing $\alpha=\sum_{i=1}^rn_i\alpha_i$ with all
$n_i$ having the same sign, we see that $\alpha_i(\xi)=0$ if
$n_i\ne0$.  So $\alpha$ is a combination of the simple roots
annihilating $\xi$.  This implies $K_\sigma=K_\xi=K_{\sigma'}$.

For $\alpha\in R_\sigma$ denote by $\alpha(\sigma)\in\{0,\pm2\pi i\}$
the constant value of $\alpha$ on $\sigma$.  Then $H_\sigma$ is the
intersection of the affine hyperplanes $H_{\alpha,\sigma(\alpha)}$
over all $\alpha\in R_\sigma$, i.e.\ 
$$
H_\sigma=\{\,\xi\in\t\mid\alpha(\xi)=\alpha(\sigma)\text{ for all
$\alpha\in R_\sigma$}\,\}.
$$
On the other hand $\z(\k_\sigma)=\bigcap\{\,\ker\alpha\mid\alpha\in
R_\sigma\,\}$, so $H_\sigma$ is a translate $\xi+\z(\k_\sigma)$ for
some $\xi\in H_\sigma$.  In fact, we can take $\xi\in\sigma$ and
therefore $\exp H_\sigma\subset Z(K_\sigma)$ since $\exp\xi\in
Z(K_\sigma)$.  If $\sigma\ge0$ then $0\in H_\sigma$, so
$H_\sigma=\z(\k_\sigma)$ and therefore $\exp H_\sigma=Z(K_\sigma)^0$.

Let $\xi$ be the closest point to the origin in the affine subspace
$H_\sigma$ (relative to the distance function given by an invariant
inner product on $\k$) and put $g_\sigma=\exp\xi$.  Then $g_\sigma\in
Z(K_\sigma)$ by \eqref{item;face-centre}.  Furthermore, $\xi$ is
perpendicular to the linear subspace parallel to $H_\sigma$, i.e.\
$\xi\in\z(\k_\sigma)^\perp=[\k_\sigma,\k_\sigma]$.  We conclude that
$g_\sigma\in[K_\sigma,K_\sigma]\cap Z(K_\sigma)$.  Since
$H_\sigma=\xi+\z(\k_\sigma)$ we have $\exp H_\sigma\subset g_\sigma
Z(K_\sigma)^0$.  The uniqueness of $g_\sigma$ up to $\Gamma_\sigma$ is
obvious.
\end{proof}

It is now a simple matter to construct a base (set of simple roots)
for each of the root systems $R_\sigma$.  Let $B$\glossary{B@$B$,
base (set of simple roots) of $R$} be the base and
$R_{\min}$\glossary{R@$R_{\min}$, minimal roots} the set of minimal
roots of $R$ determined by the chamber $\ca C$.  ($R_{\min}$ contains
one element in each irreducible component of $R$.)  By a \emph{wall}
of an open polyhedron $P$ in a vector space $V$ (e.g.\ $P$ is a
chamber or an alcove in $\t$) we mean an affine hyperplane of $V$
spanned by a codimension-one face of $P$.  Recall that the walls of
$\ca{C}\spcheck$ are the hyperplanes $H_{\alpha,0}$ with $\alpha\in
B$, while the walls of $\alcove$ include the walls of $\ca{C}\spcheck$
as well as the hyperplanes $H_{\alpha,-1}$ with $\alpha\in R_{\min}$.
Let
$$
B_\sigma=(R_\sigma\cap R_{\min})\cup\{\,\alpha\in R_\sigma\cap
B\mid\alpha(\sigma)=0\,\},
$$
where $\alpha(\sigma)$ denotes the (constant) value of any $\alpha\in
R_\sigma$ on $\sigma$.  Note that $\sigma\ge0\iff B_\sigma\subset B$.

\begin{lemma}\label{lemma;face-base}
Let $\sigma$ be a face of $\alcove$.  Then $B_\sigma$ is a base of
$R_\sigma$.  Hence $\abs{B_\sigma}=\rank[K_\sigma,K_\sigma]$ and
$\abs{B_\sigma}+\dim\sigma=\rank K$.
\end{lemma}

\begin{proof}
Let $\ca{H}_\sigma$ be the collection of singular hyperplanes of $\t$
which contain $\sigma$ and let $\ca{C}\spcheck_\sigma$ be the
connected component of $\t-\bigcup\ca{H}_\sigma$ containing $\alcove$.
Since $\exp\sigma\subset Z(K_\sigma)$, $\ca{C}\spcheck_\sigma$ is (up
to translation by any element of $\sigma$) a chamber of
$(K_\sigma,T)$.  The walls of the cone $\ca{C}\spcheck_\sigma$ are the
same as the walls of $\alcove$ which contain $\sigma$, which are
exactly the hyperplanes $H_{\alpha,\alpha(\sigma)}$ with $\alpha\in
B_\sigma$.  Therefore the elements of $B_\sigma$ are the simple roots
of $(K_\sigma,T)$ determined by the chamber $\ca{C}\spcheck_\sigma$.
In particular $B_\sigma$ spans the Lie algebra of the maximal torus
$[K_\sigma,K_\sigma]\cap T$ of $[K_\sigma,K_\sigma]$ and hence
$\abs{B_\sigma}=\rank[K_\sigma,K_\sigma]$.  By Lemma
\ref{lemma;face}\eqref{item;face-centre}, $\dim\sigma=\dim
Z(K_\sigma)$, so $\abs{B_\sigma}+\dim\sigma=\rank K_\sigma=\rank K$.
\end{proof}

Lemma \ref{lemma;face}\eqref{item;face} can be much improved for a
face of dimension greater than one.  It turns out that the centralizer
of such a face is determined by the centralizers of its proper
subfaces.  Let us write $\sigma<\tau$\glossary{.@$<$, strict inclusion
of faces} if $\sigma\le\tau$ and $\sigma\ne\tau$.

\begin{proposition}\label{proposition;centralizer}
Let $\tau$ be a face of $\alcove$ of dimension at least two.  Then
\begin{enumerate}
\item\label{item;centralizer}
$K_\tau=\bigcap_{\sigma<\tau} K_{\sigma}$;
\item\label{item;centralizer-commutator}
$[K_\tau,K_\tau]=\bigcap_{\sigma<\tau}[K_\sigma,K_\sigma]$.
\end{enumerate}
\end{proposition}

\begin{proof}
The inclusion $K_\tau\subset\bigcap_{\sigma<\tau} K_{\sigma}$ follows
from Lemma \ref{lemma;face}\eqref{item;face}.  Conversely, suppose
$g\in K_\sigma$ for all $\sigma<\tau$.  Then $\Ad(g)\exp\xi=\exp\xi$
for all $\xi$ in the (combinatorial) boundary $\partial\tau$ of
$\tau$.  Hence for each $\xi\in\partial\tau$ there exists $\gamma$ in
the exponential lattice
$\Gamma(T)=\exp_T\inv(1)$\glossary{G@$\Gamma(T)$, exponential lattice}
such that $\Ad(g)\xi=\xi+\gamma$.  The boundary of $\tau$ is connected
because $\dim\tau\ge2$.  The map $\Ad(g)$ being continuous, we
conclude that $\gamma$ is independent of $\xi\in\partial\tau$.  Let
$\eta\in\tau$ and write $\eta$ as a convex combination
$\eta=\sum_ic_i\xi_i$ with $\xi_i\in\partial\tau$ and $\sum_ic_i=1$.
Then
$$
\Ad(g)\eta=\sum_ic_i\Ad(g)\xi_i=\sum_ic_i(\xi_i+\gamma)=\eta+\gamma
$$
and therefore $\Ad(g)\exp\eta=\exp\eta$, i.e.\ $g\in K_\tau$.  Thus
$\bigcap_{\sigma<\tau} K_{\sigma}\subset K_\tau$.

Put $H=\bigcap_{\sigma<\tau}[K_\sigma,K_\sigma]$.  Then
$H\subset\bigcap_{\sigma<\tau}K_\sigma=K_\tau$ by
\eqref{item;centralizer}.  In addition Lemma
\ref{lemma;face}\eqref{item;face} implies
$[K_\tau,K_\tau]\subset[K_\sigma,K_\sigma]$ for $\sigma\le\tau$ and
hence $[K_\tau,K_\tau]\subset H$.  Now $[K_\tau,K_\tau]\subset H$
implies that $H$ (as well as $[K_\tau,K_\tau]$) is a normal subgroup
of $K_\tau$ and therefore, $T$ being a maximal torus of $K_\tau$, the
reverse inclusion $H\subset[K_\tau,K_\tau]$ will follow from
\begin{equation}\label{equation;intersection-tori}
H\cap T\subset[K_\tau,K_\tau].
\end{equation}
To prove this we may assume without loss of generality that $K$ is
almost simple, i.e.\ that $R$ is irreducible.  Let $x\in H\cap T$.
Since $\dim\tau\ge2$, $\tau$ has at least two codimension-one faces
$\rho$ and $\sigma$, so we can write $x=\exp\xi=\exp\eta$ where
$\xi\in[\k_\rho,\k_\rho]\cap\t$ and
$\eta\in[\k_\sigma,\k_\sigma]\cap\t$.  Thus $\xi-\eta$ is in the
exponential lattice $\Gamma(T)$, which is equal to the coroot lattice
$Q(R\spcheck)$\glossary{Q@$Q(R\spcheck)$, coroot lattice} since $K$ is
simply connected.  The bases of the root systems $R_\rho$ and
$R_\sigma$ are of the form $B_\rho=B_\tau\cup\{\alpha\}$ and
$B_\sigma=B_\tau\cup\{\beta\}$ where $\alpha$, $\beta\in B\cup
R_{\min}$ are two distinct roots not contained in $B_\tau$.  Hence
$\xi=\xi_0+a\alpha\spcheck$ and $\eta=\eta_0+b\beta\spcheck$ with
$\xi_0$, $\eta_0\in[\k_\tau,\k_\tau]\cap\t$ and $a$, $b\in i\R$.  Thus
\begin{equation}\label{equation;coroots}
\xi_0-\eta_0+a\alpha\spcheck-b\beta\spcheck=\xi-\eta\in Q(R\spcheck).
\end{equation}

First suppose that $\tau\ge0$.  Then we can choose $\rho$ and $\sigma$
so that $\rho\ge0$ and $\sigma\ge0$.  This implies that $\alpha$,
$\beta\in B$ are simple and therefore $a$, $b\in2\pi i\Z$ by
\eqref{equation;coroots}.  We conclude that
$x=\exp\xi_0\in[K_\tau,K_\tau]$.

Now suppose that $\tau\not\ge0$.  Then, since $R$ is irreducible,
$B_\tau$ contains the unique minimal root $\alpha_0$ and so $\alpha$
and $\beta$ are simple.  In fact, $\tau$ has at least one more
codimension-one face $\upsilon$ distinct from $\rho$ and $\sigma$; we
have $B_\upsilon=B_\tau\cup\{\gamma\}$ with $\gamma\in B$ distinct
from $\alpha$ and $\beta$; and we can write $x=\exp\zeta$ with
$\zeta=\zeta_0+c\gamma\spcheck$, $\zeta_0\in[\k_\tau,\k_\tau]\cap\t$
and $c\in i\R$.  Expanding the coroot
$-\alpha_0\spcheck=\sum_{\delta\in B}m_\delta\delta\spcheck$ in the
dual base $B\spcheck$ with positive integer coefficients $m_\delta$ we
find using \eqref{equation;coroots} that the coefficient of
$\alpha_0\spcheck$ in $\xi_0-\eta_0$ is an integer multiple of $2\pi
i/m_\gamma$ and therefore $a\in2\pi i\Z(1,m_\alpha/m_\gamma)$.
Permuting $\alpha$, $\beta$ and $\gamma$ we get the integrality
conditions
\begin{gather*}
a/2\pi i\in\Z(1,m_\alpha/m_\beta)\cap\Z(1,m_\alpha/m_\gamma),\quad
b/2\pi i\in\Z(1,m_\beta/m_\gamma)\cap\Z(1,m_\beta/m_\alpha),\\
c/2\pi i\in\Z(1,m_\gamma/m_\alpha)\cap\Z(1,m_\gamma/m_\beta).
\end{gather*}
By looking at the values of the coefficients $m_\delta$ for each
irreducible root system (see table in Example \ref{example;planches})
one can check that these conditions imply that for all triples of
distinct simple roots $\alpha$, $\beta$, $\gamma$ at least one of the
numbers $a$, $b$, $c$ is in $2\pi i\Z$ and hence that
$x\in[K_\tau,K_\tau]$.  This proves
\eqref{equation;intersection-tori}.
\end{proof}

\subsection*{Slices}
 
Using Lemma \ref{lemma;face} we can define a large ``logarithmic''
chart on $K$.  Let $O$ be the $W$-invariant open neighbourhood of the
origin in $\t$ consisting of all $\xi$ with $(2\pi
i)\inv\alpha(\xi)<1$ for all roots $\alpha$ and let
$U=\Ad(K)O=\Ad(K)\bigl(O\cap\bar{\ca{C}}\spcheck\bigr)$.

\begin{lemma}\label{lemma;chart}
\begin{enumerate}
\item\label{item;chamber-alcove}
$O\cap\ca{C}\spcheck=\alcove$ and $O\cap\bar{\ca{C}}\spcheck
=\bigcup\{\,\sigma\mid0\le\sigma\le\alcove\,\}$.
\item\label{item;etale-chart}
$\exp U$ is open and dense in $K$ and $\exp\colon U\to\exp U$ is a
local diffeomorphism.
\item\label{item;chart}
$\exp\colon\frac12U\to K$ is a diffeomorphism onto its image.
\end{enumerate}
\end{lemma}

\begin{proof}
Let $R_+$\glossary{R@$R_+$, positive roots} denote the set of positive
roots with respect to the chamber $\ca{C}$.  By \cite[Ch.~6, \S~2.3,
Proposition 5]{bourbaki;groupes-algebres} a point $\xi\in\t$ is in the
alcove $\alcove$ if and only if $0<(2\pi i)\inv\alpha(\xi)<1$ for all
$\alpha\in R_+$.  Let us abbreviate this system of inequalities for
$\alcove$ to
\begin{equation}\label{equation-alcove}
0<(2\pi i)\inv\alpha<1\quad\text{for $\alpha\in R_+$}.
\end{equation}
Likewise, the chamber $\ca{C}\spcheck$ in $\t$ is given by $(2\pi
i)\inv\alpha>0$ for $\alpha\in R_+$.  Thus $\xi\in
O\cap\ca{C}\spcheck$ iff $0<(2\pi i)\inv\alpha(\xi)<1$ for all
$\alpha\in R_+$ iff $\xi\in\alcove$.  Hence
$O\cap\ca{C}\spcheck=\alcove$.  This implies
$O\cap\bar{\ca{C}}\spcheck\subset\bar\alcove$; in fact
$$
\xi\in O\cap\bar{\ca{C}}\spcheck\iff0\le(2\pi
i)\inv\alpha(\xi)<1\quad\text{for all $\alpha\in R_+$}.
$$
So if $\sigma\le\alcove$ is the face containing $\xi\in
O\cap\bar{\ca{C}}\spcheck$ then $\alpha(\sigma)=0$ for all $\alpha\in
R_\sigma$, and hence $\sigma\ge0$.  Conversely, if
$\xi\in\bar{\alcove}$ is contained in a face $\sigma\ge0$, then $(2\pi
i)\inv\alpha(\xi)<1$ for all $\alpha\in R_+$, so $\xi\in
O\cap\bar{\ca{C}}\spcheck$.  Thus $O\cap\bar{\ca{C}}\spcheck
=\bigcup\{\,\sigma\mid0\le\sigma\le\alcove\,\}$.

The formula $T_\xi\exp=T_1L(\exp\xi)\circ(1-e^{-\ad\xi})/\ad\xi$ shows
that $T_\xi\exp$ is bijective if and only if $\ad\xi$ does not have an
eigenvalue $2\pi in$ with $n\in\Z-\{0\}$.  This implies that
$T_\xi\exp$ is bijective for every point in $U$.  Therefore $\exp U$
is open and $\exp\colon U\to\exp U$ is a local diffeomorphism.  Since
$\exp U$ contains the set of regular elements $\Ad(K)\exp\alcove$, it
is dense in $K$.

It remains to show that $\exp\colon\frac12U\to K$ is injective.
Suppose $\xi$, $\eta\in\frac12U$ satisfy $\exp\xi=\exp\eta$.  Without
loss of generality we may assume
$\xi\in\frac12O\cap\bar{\ca{C}}\spcheck$ and $\eta=\Ad(g)\xi'$ with
$g\in K$ and $\xi'\in\frac12O\cap\bar{\ca{C}}\spcheck$.  Then
$\exp\xi$, $\exp\xi'\in T$, so $\exp\xi=\Ad(g)\exp\xi'$ implies that
$\exp\xi=\Ad(n)\exp\xi'$ for some $n\in N(T)$.  Then $gn\inv$ is in
the centralizer of $\exp\xi'$, which by \eqref{item;chamber-alcove}
and Lemma \ref{lemma;face}\eqref{item;face-face} is equal to the
centralizer of $\xi'$.  Therefore
\begin{equation}\label{equation;gn}
\eta=\Ad(g)\xi'=w\xi',
\end{equation}
where $w$ is the Weyl group element represented by $n$.  From
$\exp\xi=\Ad(n)\exp\xi'$ we also get $\xi=w\xi'+\gamma$ for some
$\gamma\in\Gamma(T)$, where $\Gamma(T)=\exp_T\inv(1)$ is the
exponential lattice.  Since $\xi$, $\xi'\in\frac12O\cap\bar{\ca{C}}$
we have
$$
-1<\frac{\alpha(w\xi')}{\pi i}<1,\quad-1<\frac{\alpha(\xi)}{\pi i}
=\frac{\alpha(w\xi')+\alpha(\gamma)}{\pi i}<1
$$
for all $\alpha\in R_+$.  Since $\alpha(\gamma)\in2\pi i\Z$ this
implies $\alpha(\gamma)=0$ for all $\alpha$, so $\gamma=0$ and
$\xi=w\xi'$.  Together with \eqref{equation;gn} this shows $\xi=\eta$.
\end{proof}

Let
$\star\sigma=\bigcup_{\tau\ge\sigma}\tau$\glossary{star@$\star\sigma$,
star of face $\sigma$} denote the open star of a face
$\sigma\le\alcove$ and put
$$
\eu{S}_\sigma=\Ad(K_\sigma)\exp(\star\sigma)\subset
K_\sigma.\glossary{S@$\eu{S}_\sigma$, slice at face $\sigma$}
$$
(Note that $\star\sigma=O\cap\bar{\ca{C}}\spcheck$ and
$\eu{S}_\sigma=U$ for $\sigma=0$.)  The next result (stated without
proof in \cite{alekseev-malkin-meinrenken;lie-group}) says that
$\eu{S}_\sigma$ is a slice for the conjugation action of $K$ at all
points in $\exp\sigma$.

\begin{proposition}\label{proposition;slice}
Let $\sigma$ be a face of $\alcove$.  Then $\eu{S}_\sigma$ is a
submanifold of $K$.  Its orbit $\Ad(K)\eu{S}_\sigma$ under the
conjugation action is open in $K$.  The conjugation map
$K\times\eu{S}_\sigma\to K$ descends to a diffeomorphism from the
associated bundle $K\times^{K_\sigma}\eu{S}_\sigma$ onto
$\Ad(K)\eu{S}_\sigma$.  In particular, if $g\in\eu{S}_\sigma$ and
$kgk\inv\in\eu{S}_\sigma$, then $k\in K_\sigma$.
\end{proposition}

\begin{proof}
The face $\sigma$ and its star can be described by inequalities
similar to \eqref{equation-alcove}.  Let $R_{\sigma,+}=R_\sigma\cap
R_+$ be the set of positive roots in $R_\sigma$.  Then $\sigma$ is
given by
\begin{equation}\label{equation;inequality-face}
\begin{cases}
0<(2\pi i)\inv\alpha<1&\text{for $\alpha\in R_+-R_{\sigma,+}$},\\
\alpha=\alpha(\sigma)&\text{for $\alpha\in R_{\sigma,+}$},
\end{cases}
\end{equation}
while $\star\sigma$ is given by
\begin{equation}\label{equation;inequality-star}
\begin{cases}
0<(2\pi i)\inv\alpha<1&\text{for $\alpha\in R_+-R_{\sigma,+}$},\\
0\le(2\pi i)\inv\alpha<1&\text{for $\alpha\in R_{\sigma,+}$,
$\alpha(\sigma)=0$},\\
0<(2\pi i)\inv\alpha\le1&\text{for $\alpha\in R_{\sigma,+}$,
$\alpha(\sigma)=2\pi i$}.\\
\end{cases}
\end{equation}
Let $\alcove_\sigma\subset\t$ be the unique alcove for the subgroup
$K_\sigma$ containing $\alcove$.  By \eqref{equation-alcove},
$\alcove_\sigma$ is given by the inequalities $0<(2\pi i)\inv\alpha<1$
for $\alpha\in R_{\sigma,+}$.  Hence its closure $\bar{\alcove}_\sigma$
is given by $0\le(2\pi i)\inv\alpha\le1$ for $\alpha\in R_{\sigma,+}$.
Comparing with \eqref{equation;inequality-star} we see that
$\star\sigma$ is open in $\bar{\alcove}_\sigma$.  Hence the image of
$\star\sigma$ in the quotient $\bar{\alcove}_\sigma/\pi_1(K_\sigma)$ is
open.  By \cite[Ch.~9, \S~5.2, Corollaire
1]{bourbaki;groupes-algebres}, the exponential map
$\bar{\alcove}_\sigma\to T$ induces a homeomorphism from
$\bar{\alcove}_\sigma/\pi_1(K_\sigma)$ to $K_\sigma/\Ad(K_\sigma)$.  We
conclude that $\eu{S}_\sigma$ is an open subset of $K_\sigma$, and
therefore an embedded submanifold of $K$.  Let $\mu\colon
K\times\eu{S}_\sigma\to K$ be the conjugation map and let
$g\in\eu{S}_\sigma$.  A tangent vector in
$T_{(1,g)}(K\times\eu{S}_\sigma)=\k\oplus T_gK_\sigma$ is of the form
$(\xi,R(g)_*\eta)$ with $\xi\in\k$ and $\eta\in\k_\sigma$.  A
computation yields
$$
T_{(1,g)}\mu(\xi,R(g)_*\eta)=R(g)_*\xi-L(g)_*\xi+R(g)_*\eta.
$$
Hence the kernel of $T_{(1,g)}\mu$ consists of all $(\xi,R(g)_*\eta)$
with
\begin{equation}\label{equation;kernel}
\eta=\Ad(g)\xi-\xi.
\end{equation}
Since $g\in\eu{S}_\sigma$, we can write $g=\exp(\Ad(k)\zeta)$ with
$k\in K_\sigma$ and $\zeta\in\star\sigma$.  Then
\eqref{equation;kernel} is equivalent to
$\eta'=e^{\ad\zeta}\xi'-\xi'$, where $\eta'=\Ad(k)\inv\eta$ and
$\xi'=\Ad(k)\inv\xi$.  Decomposing $\xi'=h+\sum_{\alpha\in R}c_\alpha
x_\alpha$ with $h\in\t$ we find
\begin{equation}\label{equation;kalpha}
\eta'=e^{\ad\zeta}\xi'-\xi'
=\sum_{\alpha}c_\alpha\bigl(e^{\alpha(\zeta)}-1\bigr)x_\alpha.
\end{equation}
Since $\zeta\in\star\sigma$, we have $0<(2\pi i)\inv\alpha(\zeta)<1$
for all $\alpha\in R_+-R_{\sigma,+}$ by
\eqref{equation;inequality-star}.  Also $\eta'\in\k_\sigma$, so by
\eqref{equation;kalpha} we have $c_\alpha=0$ for all $\alpha\in
R-R_\sigma$.  This means $\xi'\in\k_\sigma$, and therefore
$\xi=\Ad(k)\xi'\in\k_\sigma$.  Putting this together with
\eqref{equation;kernel} we see that
$$
\ker
T_{(1,g)}\mu=\{\,(\xi,L(g)_*\xi-R(g)_*\xi)\mid\xi\in\k_\sigma\,\},
$$
which is exactly the kernel of the $K_\sigma$-quotient map
$K\times\eu{S}_\sigma\to K\times^{K_\sigma}\eu{S}_\sigma$.  This shows
that the map $\bar{\mu}\colon K\times^{K_\sigma}\eu{S}_\sigma\to K$
induced by $\mu$ is immersive at all points of the form $(1,g)$.
Being $K$-equivariant $\bar{\mu}$ is an immersion at every point.
Moreover, $\dim K\times^{K_\sigma}\eu{S}_\sigma=\dim K+\dim
K_\sigma-\dim K_\sigma=\dim K$, so $\Ad(K)\eu{S}_\sigma$ is open in
$K$ and $\bar{\mu}$ is a local diffeomorphism.  It remains to show
that $\bar{\mu}$ is injective, which is equivalent to the last
statement of the proposition.  So let $k\in K$, $g\in\eu{S}_\sigma$
and suppose $kgk\inv\in\eu{S}_\sigma$.  Without loss of generality we
may assume $g=\exp\xi$ with $\xi\in\star\sigma$.  We can write
$k(\exp\xi)k\inv=h(\exp\eta)h\inv$ for certain $h\in K_\sigma$ and
$\eta\in\star\sigma$.  Both $\xi$ and $\eta$ are in $\bar{\alcove}$, so
from the identification $\bar{\alcove}\cong K/\Ad(K)$ we obtain
$\eta=\xi$, whence $h\inv k\in K_{\exp\xi}$.  Since
$K_{\exp\xi}\subset K_\sigma$ by Lemma
\ref{lemma;face}\eqref{item;face-face} and $h\in K_\sigma$, we
conclude $k\in K_\sigma$.
\end{proof}

\subsection*{The action of the centre}
 
The translation action of the centre $Z(K)$ on $K$ commutes with the
adjoint action of $K$ and so descends to an action on the space of
conjugacy classes.  Via the identification $K/\Ad
K\cong\exp\bar{\alcove}$ we can think of this action as taking place
on $\exp\bar{\alcove}$, namely by shifting $\exp\bar{\alcove}$ by a
central element and reflecting the resulting subset of $T$ back into
$\exp\bar{\alcove}$ by a Weyl group element.

\begin{lemma}\label{lemma;centre}
For each $c\in Z(K)$ there exists a unique $w$ in the Weyl group $W$
such that $w(L(c)\exp\alcove)=\exp\alcove$.  The map $\zeta\colon
Z(K)\to W$ sending $c$ to $w$ is an injective homomorphism.
\end{lemma}

\begin{proof}
By \cite[Ch.~9, \S~5.4, Proposition 5]{bourbaki;groupes-algebres} the
inverse image of the centre $Z(K)$ under $\exp_T\colon\t\to T$ is the
coweight lattice $P(R\spcheck)\subset\t$\glossary{P@$P(R\spcheck)$,
coweight lattice}, the lattice dual to the root lattice in $\t^*$
(whether or not $K$ is simply connected).  In other words
\begin{equation}\label{equation;centre}
Z(K)=\{\,\exp\xi\mid\text{$\xi\in\t$ and $\alpha(\xi)\in2\pi i\Z$ for
all $\alpha\in R$}\,\}.
\end{equation}
Accordingly, let $c$ be central and write $c=\exp\xi$ with $\xi\in
P(R\spcheck)$.  If $H_{\alpha,n}=\{\,\eta\in\t\mid\alpha(\eta)=2\pi
in\,\}$\glossary{H@$H_{\alpha,n}$, singular hyperplane of $\t$} is a
singular hyperplane of $\t$, then
$\xi+H_{\alpha,n}=H_{\alpha,n+\alpha(\xi)}$, so $P(R\spcheck)$
preserves the collection of singular hyperplanes.  Therefore
$P(R\spcheck)$ maps alcoves to alcoves, so there is an element
$\tilde{w}$ of the affine Weyl group $W\affine$ mapping $\xi+\alcove$
to $\alcove$.  Since $K$ is simply connected we have
$W\affine=W\ltimes\Gamma(T)$.  Hence we can write
$\tilde{w}=wt_\gamma$ with $w\in W$ and $\gamma\in\Gamma(T)$, where
$t_\gamma$ denotes the translation by $\gamma$.  From
$\tilde{w}(\xi+\alcove)=\alcove$ we obtain
$$
w\bigl(L(c)\exp\alcove\bigr)=w\bigl(\exp(\xi+\gamma+\alcove)\bigr)
=\exp\bigl(\tilde{w}(\gamma+\alcove)\bigr)=\exp\alcove.
$$
Since $\exp\alcove$ and $L(c)\exp\alcove$ are connected components of
the set of regular points of $T$ and $W$ acts simply transitively on
those components, $w$ is unique and we can define $\zeta\colon Z(K)\to
W$ by $\zeta(c)=w$.  The left $Z(K)$-action on $T$ commutes with the
$W$-action, so from $w\bigl(L(c)\exp\alcove\bigr)=\exp\alcove$ we get
$w(\exp\alcove)=L(c)\inv\exp\alcove$.  Hence if $c'$ is another
central element,
\begin{multline*}
\zeta(cc')(\exp\alcove)=L(cc')\inv\exp\alcove=L(c')\inv
L(c)\inv\exp\alcove=L(c')\inv\zeta(c)(\exp\alcove)\\
=\zeta(c)\bigl(L(c')\inv\exp\alcove\bigr)
=\zeta(c)\zeta(c')(\exp\alcove),
\end{multline*}
which shows that $\zeta$ is a homomorphism.  If $w=\zeta(c)=1$ then
$\tilde{w}=wt_\gamma=t_\gamma$ is a translation and hence
$\gamma+\xi+\alcove=\alcove$.  This implies $\xi=-\gamma$, so
$\xi\in\Gamma(T)$ and $c=\exp\xi=1$.  Therefore $\zeta$ is injective.
\end{proof}

\begin{example}\label{example;centre-weyl}
Let $K=\SU(l+1)$ and let $T\cong\U(1)^l$ be the standard maximal
torus.  Identify the Cartan subalgebra $\t$ with the subspace
$V=\{\,x\mid x_1+x_2+\cdots+x_{l+1}=0\,\}$ of $\R^{l+1}$ via the
isomorphism $x\mapsto2\pi i\diag(x_1,x_2,\dots,x_{l+1})$.  The simple
roots are then $(2\pi i)\inv\alpha_j(x)=x_j-x_{j+1}$ for $1\le j<l+1$,
the minimal root is $(2\pi i)\inv\alpha_0(x)=-x_1+x_{l+1}$, so the
alcove is the open simplex
$$
\alcove=\{\,x\in V\mid x_{l+1}+1>x_1>x_2>\cdots>x_{l+1}\,\}.
$$
Let $\eps_1$, $\eps_2$,\dots, $\eps_{l+1}$ be the standard basis of
$\R^{l+1}$.  Then the vertices of $\alcove$ are
$$
\sigma_j=\sum_{k=1}^j\eps_k-\frac{j}{l+1}\sum_{k=1}^{l+1}\eps_k
\qquad\text{for $j=0$, $1$,\dots, $l$.}
$$
(Note that under the isomorphism $\R^{l+1}\to(\R^{l+1})^*$ given by
the standard inner product the vertex $\sigma_j$ is mapped to the
fundamental weight $\varpi_j$.)  The centre of $K$ is the cyclic group
generated by $c=\exp\sigma_1=e^{-2\pi i/(l+1)}I$; we have
$c^j=\exp\sigma_j$.  The barycentre of $\alcove$ is
$$
\beta=\frac1{l+1}\sum_{j=0}^l\sigma_j
=\frac1{2l+2}(l,l-2,\dots,2-l,-l).
$$
It is straightforward to check that $w(\gamma+\sigma_1+\beta)=\beta$,
where $\gamma=-\eps_1+\eps_{l+1}\in\Gamma(T)$ and $w$ is the cyclic
permutation $w(x)=(x_{l+1},x_1,x_2,\dots,x_l)$.  Thus $\zeta(c)=w$.
\end{example}

\section{The spinning $2n$-sphere}\label{section;sphere}

Here we describe a quasi-Hamiltonian $\U(n)$-structure on $S^{2n}$,
which generalizes the $\U(2)$-structure on $S^4$ of \cite[Appendix
A]{alekseev-meinrenken-woodward;duistermaat-heckman} and
\cite[Proposition 2.29]{hurtubise-jeffrey;representations-weighted}.
We discovered this structure in the course of our work on the imploded
cross-section of $D\U(n)$.  Because the final result is easy to
explain without employing the machinery of implosion we present it in
this appendix.  Its relationship to implosion is clarified in Example
\ref{example;sun}.

We will obtain $S^{2n}$ by gluing together two copies of a
quasi-Hamiltonian $2n$-disc.  This is similar to the construction of
$S^4$ in \cite{alekseev-meinrenken-woodward;duistermaat-heckman}, but
the generalization from $S^4$ to $S^{2n}$ is not entirely
straightforward.  To begin with, the gluing formula in [\emph{loc.\
cit.}] does not generalize to higher dimensions and must be adapted.
Moreover, it turns out that there is a second generalization of the
$\U(2)$-structure on $S^4$, namely to a
$\U(n,\H)\times\U(1)$-structure on the quaternionic projective space
$\H\P^n$.  This will be the subject of a future publication.

Let $T$ be the standard maximal torus of $\U(n)$.  As invariant inner
product on $\u(n)$ we take $(\xi,\eta)=-(4\pi^2)\inv\tr(\xi\eta)$.
Consider $\C^n$ with its standard Hermitian structure
$z^*w=\sum_j\bar{z}_jw_j$ and its standard $\U(n)$-action.  The
corresponding symplectic form $\omega_0$ and moment map $\Phi_0$ are
given by
$$
\omega_0=\Im\langle\cdot,\cdot\rangle=\frac{i}2\sum_jdz_j\,d\bar{z}_j,
\quad \Phi_0(z)=-2\pi^2i\,zz^*,
$$
where we have used the inner product on $\u(n)$ to write $\Phi_0$ as a
map into $\u(n)$.  To produce a quasi-Hamiltonian structure we
exponentiate: $\Phi=\exp\circ\Phi_0$ and
$\omega=\omega_0+\Phi_0^*\varpi$.  We determine where $\omega$ is
minimally degenerate and write a more explicit formula for it as
follows.

\begin{lemma}\label{lemma;disc}
Let $D$ be the disc $\{\,z\in\C^n\mid\norm{z}^2<1/\pi\,\}$.
\begin{enumerate}
\item\label{item;regular}
$T_\xi\exp$ is surjective for all $\xi\in\Phi_0(D)$.
\item\label{item;disc}
$(D,\omega,\Phi)$ is a quasi-Hamiltonian $\U(n)$-manifold.
\item\label{item;forms}
Let $\lambda$ be the $\U(n)$-invariant two-form on $\C^n-\{0\}$ given
by
$$
\lambda=d\log\norm{z}^2\wedge\frac12\Im(z^*\,dz)
=\frac{i}{2\norm{z}^2}\sum_{1\le i\le k\le
n}\bar{z}_jz_k\,dz_j\,d\bar{z}_k.
$$
Then $\displaystyle\omega=\lambda
+\frac{\sin(2\pi^2\norm{z}^2)}{2\pi^2\norm{z}^2}(\omega_0-\lambda)$.
\end{enumerate}
\end{lemma}

Note that $\omega$ is smooth even though $\lambda$ has a pole at the
origin.  In fact, the value of $\omega$ at the origin is equal to the
standard symplectic form $\omega_0$.

\begin{proof}
The Hermitian matrix $zz^*$ has spectrum $(\norm{z}^2,0,\dots,0)$.
Hence the matrix $\Phi_0(z)$ is conjugate to the diagonal matrix
$\xi=-2\pi^2i\diag(\norm{z}^2,0,\dots,0)\in\t$.  Thus if $\alpha$ is a
root we have either $\alpha(\xi)=\pm2\pi^2i\norm{z}^2$ or
$\alpha(\xi)=0$.  According to Lemma
\ref{lemma;chart}\eqref{item;etale-chart} $T_\xi\exp$ is surjective if
$(2\pi i)\inv\alpha(\xi)<1$ for all $\alpha$, i.e.\ if $z\in D$.  This
proves \eqref{item;regular}.  It follows that $\omega$ is minimally
degenerate on $D$, which proves \eqref{item;disc}.  Because of
$\U(n)$-invariance, for the proof of \eqref{item;forms} it suffices to
consider vectors $z$ of the form $(z_1,0,\dots,0)$.  At such a point
$z$ we have
\begin{equation}\label{equation;lambda}
\lambda=\frac{i}2\,dz_1\,d\bar{z}_1
\end{equation}
and so
\begin{equation}\label{equation;lambda-omega}
\lambda+\frac{\sin(2\pi^2\norm{z}^2)}{2\pi^2\norm{z}^2}
(\omega_0-\lambda)=\frac{i}2\,dz_1\,d\bar{z}_1
+\frac{i\sin(2\pi^2\norm{z}^2)}{4\pi^2\norm{z}^2}
\sum_{r=2}^ndz_r\,d\bar{z}_r.
\end{equation}
To compute $\omega$ we first calculate $\varpi\in\Omega^2(\u(n))$, or
rather its complex bilinear extension to an element of
$\Omega_\C^2(\lie{gl}(n,\C))$.  Let $\lambda\in\t$.  Let
$\{\eps_{jk}\}_{1\le j,k\le n}$ be the standard basis of
$\lie{gl}(n,\C)$, with dual basis $\{d\xi_{jk}\}_{1\le j,k\le n}$ in
$\lie{gl}(n,\C)^*$, and for $j\ne k$ let $\alpha_{jk}$ be the root
$d\xi_{jj}-d\xi_{kk}$.  It follows from \eqref{equation;varpi} that
$\varpi_\lambda(\eps_{jj},\eps_{kk})=0$ for all $j$ and $k$.  Likewise
for $j\ne k$ and $l\ne m$ we have
$$
\varpi_\lambda(\eps_{jk},\eps_{lm})
=\bigl(f(\alpha_{jk}(\lambda))\eps_{jk},\eps_{lm}\bigr)
=-\frac{f(\alpha_{jk}(\lambda))}{4\pi^2}\delta_{kl}\delta_{jm},
$$
where $f(x)=(x-\sinh x)/x^2$.  Using
$$
f(\alpha_{jk}(\lambda))\delta_{kl}\delta_{jm}
=\sum_{r<s}f(\alpha_{rs}(\lambda))\,
d\xi_{rs}\,d\xi_{sr}(\eps_{jk},\eps_{lm})
$$
we conclude that for all $\lambda\in\t$
\begin{equation}\label{equation;varpi-lambda}
\varpi_\lambda
=-\frac1{4\pi^2}\sum_{r<s}f(\alpha_{rs}(\lambda))\,d\xi_{rs}\,d\xi_{sr}.
\end{equation}
The pullback $(\Phi_0^*\varpi)_z$ is found by substituting
$\lambda=\Phi_0(z)$ and $\xi_{rs}=-2\pi^2iz_r\bar{z}_s$ into $\varpi$.
This yields $\alpha_{rs}(\lambda)=-2\pi^2i(\abs{z_r}^2-\abs{z_s}^2)$
and
$$
d\xi_{rs}\,d\xi_{sr}=-4\pi^4\bigl(\abs{z_s}^2\,dz_r\,d\bar{z}_r
-\abs{z_r}^2\,dz_s\,d\bar{z}_s
+2i\Im(\bar{z}_r\bar{z}_s\,dz_r\,dz_s)\bigr),
$$
whence 
\begin{multline*}
(\Phi_0^*\varpi)_z=\pi^2\sum_{1\le r<s\le
n}f\bigl(-2\pi^2i\bigl(\abs{z_r}^2-\abs{z_s}^2\bigr)\bigr)
\bigl(\abs{z_s}^2\,dz_r\,d\bar{z}_r\\-\abs{z_r}^2\,dz_s\,d\bar{z}_s
+2i\Im(\bar{z}_r\bar{z}_s\,dz_r\,dz_s)\bigr).
\end{multline*}
Since $z_2=\cdots=z_n=0$ this simplifies to
$$
(\Phi_0^*\varpi)_z=-\pi^2\norm{z}^2f\bigl(-2\pi^2i\norm{z}^2\bigr)
\sum_{s=2}^ndz_s\,d\bar{z}_s.
$$
From $f(x)=(x-\sinh x)/x^2$ we get
$$
(\Phi_0^*\varpi)_z=\biggl(-\frac{i}2
+\frac{i\sin(2\pi^2\norm{z}^2)}{4\pi^2\norm{z}^2}\biggr)
\sum_{s=2}^ndz_s\,d\bar{z}_s
$$
and so
$$
\omega_z=(\omega_0)_z+(\Phi_0^*\varpi)_z=\frac{i}2\,dz_1\,d\bar{z}_1
+\frac{i\sin(2\pi^2\norm{z}^2)}{4\pi^2\norm{z}^2}
\sum_{r=2}^ndz_r\,d\bar{z}_r.
$$
Comparing with \eqref{equation;lambda-omega} we obtain
\eqref{item;forms}.
\end{proof}

Now define a $\U(n)$-equivariant diffeomorphism from the punctured
disc $D-\{0\}$ to itself by
\begin{equation}\label{equation;transition}
\phi(z)=-s(z)z,\qquad\text{where}\quad
s(z)=\frac{\sqrt{\pi\inv-\norm{z}^2}}{\norm{z}}.
\end{equation}

\begin{lemma}\label{lemma;glue-disc}
$\phi^*\omega=-\omega$ and $\phi^*\Phi=i\circ\Phi$.
\end{lemma}

\begin{proof}
Let $p(z)$ be the anti-Hermitian matrix $2\pi izz^*/\norm{z}^2$.  Then
\begin{equation}\label{equation;pull}
\Phi_0(\phi(z))=-2\pi^2i\,s(z)^2zz^*=-p(z)-\Phi_0(z).
\end{equation}
Note that $p(z)$ is $2\pi i$ times the orthogonal projection onto the
line spanned by $z$.  Therefore $\exp p(z)=1$, whence
$\Phi(\phi(z))=\exp(\Phi_0(\phi(z)))=\exp(-\Phi_0(z))=\Phi(z)\inv$.
From \eqref{equation;pull} we also see
$$
\phi^*\omega
=\phi^*(\omega_0+\Phi_0^*\varpi)=\phi^*\omega_0-p^*\varpi-\Phi_0^*\varpi.
$$
A straightforward computation yields
$\phi^*\omega_0=s(z)^2\omega_0-\pi\inv\norm{z}^{-2}\lambda$, and so
\begin{align*}
\phi^*\omega&
=s(z)^2\omega_0-\frac1{\pi\norm{z}^2}\lambda-p^*\varpi-\Phi_0^*\varpi\\
&=-(\omega_0+\Phi_0^*\varpi)+\frac1{\pi\norm{z}^2}(\omega_0-\lambda)
-p^*\varpi\\
&=-\omega+\frac1{\pi\norm{z}^2}(\omega_0-\lambda)
-p^*\varpi.
\end{align*}
It remains to show that
\begin{equation}\label{equation;p-varpi}
p^*\varpi=\frac1{\pi\norm{z}^2}(\omega_0-\lambda).
\end{equation}
Because both sides are $\U(n)$-invariant forms, it suffices to check
this identity at points $z\in\C^n$ of the form $(z_1,0,\dots,0)$.
From \eqref{equation;lambda} we see
\begin{equation}\label{equation;omega-lambda}
\frac1{\pi\norm{z}^2}(\omega_0-\lambda)
=\frac{i}{2\pi\norm{z}^2}\sum_{r=2}^ndz_r\,d\bar{z}_r.
\end{equation}
To find $(p^*\varpi)_z$ we substitute $\lambda=\Phi_0(z)$ and
$\xi_{rs}=2\pi iz_r\bar{z}_s/\norm{z}^2$ into
\eqref{equation;varpi-lambda}, which yields
$$
(p^*\varpi)_z=-\sum_{r=2}^nf(2\pi i)\norm{z}^{-2}\,dz_r\,d\bar{z}_r=
\frac{i}{2\pi\norm{z}^2}\sum_{r=2}^ndz_r\,d\bar{z}_r.
$$
Together with \eqref{equation;omega-lambda} this gives
\eqref{equation;p-varpi}.
\end{proof}

This lemma tells us that the quasi-Hamiltonian disc $(D,\omega,\Phi)$
and its opposite $(D,-\omega,i\circ\Phi)$ can be glued together
smoothly via the map $\phi$, resulting in a quasi-Hamiltonian
structure on the $2n$-sphere.  Because the gluing map $\phi$ sends $D$
to itself, it can also be viewed as a diffeomorphism from $S^{2n}$ to
itself, which reverses the quasi-Hamiltonian structure.  Under a
suitable identification of $D\amalg_\phi D$ with $S^{2n}$ this map is
nothing but the antipodal map.  To summarize,

\begin{theorem}\label{theorem;glue-disc}
The $2n$-dimensional sphere with its standard $\U(n)$-action is a
quasi-Hamiltonian $\U(n)$-manifold.  The antipodal map is an
anti-automorphism of the quasi-Hamiltonian structure.
\end{theorem}

Define the maps $\Phi_1\colon D\to\SU(n)$ and $\Phi_2\colon D\to
Z\U(n)$ by
$$
\Phi_1(z)=\exp2\pi^2i\biggl(\frac{\norm{z}^2}{n}I-zz^*\biggr),
\qquad\Phi_2(z)=e^{-2\pi^2i\norm{z}^2/n}I
$$
and observe that $\Phi=\Phi_1\Phi_2$.  Hence, by Lemma
\ref{lemma;quasi}\eqref{item;cover}, $S^{2n}$ is a quasi-Hamilton\-ian
$\SU(n)\times Z\U(n)$-manifold with moment map $\Phi_1\times\Phi_2$.
By Theorem \ref{theorem;quotient}, the reduction of $S^{2n}$ with
respect to the central circle $Z\U(n)$ is a quasi-Hamiltonian
$\SU(n)$-manifold.  Reduction at the poles $z=0$ and
$\norm{z}^2=\pi\inv$ gives a quotient consisting of a single point,
while reduction at an intermediate level $\norm{z}^2=a$ gives a
projective space $\C\P^{n-1}$.  Note that the form $\lambda$ of Lemma
\ref{lemma;disc}\eqref{item;forms} vanishes on the $2n-1$-sphere
$\norm{z}^2=a$.  Therefore the induced two-form on the quotient is
$$
\frac{\sin(2\pi^2a)}{2\pi^2a}\sigma,
$$
where $\sigma$ is a multiple of the Fubini-Study symplectic form.
Thus we have found a family of quasi-Hamiltonian $\SU(n)$-structures
on $\C\P^{n-1}$ for which the two-forms happen to be closed and
non-degenerate (except for $a=1/2\pi$, which corresponds to the
equator of $S^{2n}$, where the induced form is zero).  It is easy to
see (e.g.\ by using Addendum \ref{addendum;quotient}) that these
quasi-Hamiltonian $\SU(n)$-structures are the same as those one
obtains by considering $\C\P^{n-1}$ as a conjugacy class of $\SU(n)$.

\newpage


 \begin{theglossary}\begin{multicols}{2}

  \item $<$, strict inclusion of faces, 30
  \item $\le$, inclusion of faces, 6, 8
  \item $\fusion$, fusion product, 4
  \item $\sim$, equivalence relation defining $M\impl$, 7

  \indexspace

  \item $\alcove$, open alcove in $\t$, 7
  \item $\bar{\alcove}$, closure of $\alcove$, 7
  \item $\eu{A}$, $\hat{T}$-action on $DK$, 15
  \item $A^0$, unit component of a Lie group $A$, 29
  \item $\Ad$, adjoint action on $K$ or $\k$, 3

  \indexspace

  \item $B$, base (set of simple roots) of $R$, 29
  \item $\eu{B}$, $T_0$-action on $DK$, 16

  \indexspace

  \item $\ca{C}$, open chamber in $\t^*$, 6
  \item $\bar{\ca{C}}$, closure of $\ca{C}$, 6
  \item $\ca{C}\spcheck$, dual open chamber in $\t$, 7
  \item $\chi$, bi-invariant three-form, 3

  \indexspace

  \item $DK$, double $K\times K$, 12

  \indexspace

  \item $G$, complexification of $K$, 21
  \item $\Gamma(T)$, exponential lattice, 30
  \item $\Gamma_\sigma$, covering group $[K_\sigma,K_\sigma]\cap Z(K_\sigma)^0$, 
		9, 29

  \indexspace

  \item $H_{\alpha,n}$, singular hyperplane of $\t$, 28, 33

  \indexspace

  \item i, inversion $g\mapsto g\inv$, 4
  \item impl, implosion, 7

  \indexspace

  \item $K$, compact simply connected Lie group, 3, 7
  \item $K_\sigma$, centralizer of a face, 7, 29
  \item $K_g$, centralizer of $g\in K$, 7
  \item $\k_\sigma$, Lie algebra of $K_\sigma$, 29

  \indexspace

  \item $L(g)$, left multiplication by $g$, 3

  \indexspace

  \item $M$, (quasi-)Hamiltonian $K$-manifold, 3
  \item $M\impl$, imploded cross-section, 7
  \item $M\quot[g]K$, (quasi-)Hamiltonian quotient, 5
  \item $M_\sigma$, cross-section over $\sigma$, 7, 8

  \indexspace

  \item $N$, maximal unipotent subgroup of $G$, 21

  \indexspace

  \item $\omega$, (quasi-)symplectic form, 3

  \indexspace

  \item $P(R\spcheck)$, coweight lattice, 33
  \item $\Phi$, moment map, 3
  \item $\Phi\impl$, imploded moment map, 7
  \item $\Phi_\sigma$, moment map on $M_\sigma$, 8

  \indexspace

  \item $Q(R\spcheck)$, coroot lattice, 30

  \indexspace

  \item $R$, root system of $(K,T)$, 28
  \item $R_+$, positive roots, 31
  \item $R_\sigma$, root system of $(K_\sigma,T)$, 28
  \item $R_{\min}$, minimal roots, 29
  \item $R(g)$, right multiplication by $g$, 3

  \indexspace

  \item $\eu{S}_\sigma$, slice at face $\sigma$, 7, 8, 32
  \item $\sigma$, face of chamber or alcove, 6, 28
  \item $\sigma\prin$, principal face, 6, 8
  \item $\star\sigma$, star of face $\sigma$, 7, 8, 32

  \indexspace

  \item $T$, maximal torus of $K$, 6
  \item $T_0$, extension of $\Z/2\Z$ by maximal torus, 16
  \item $\hat{T}$, extension of centre by maximal torus, 14
  \item $\theta_L$, $\theta_R$, Maurer-Cartan forms, 3

  \indexspace

  \item $W$, Weyl group, 7
  \item $w_0$, longest Weyl group element, 16

  \indexspace

  \item $X_\sigma$, piece $\Phi\inv(\exp\sigma)\big/[K_\sigma,K_\sigma]$ of $M\impl$, 
		8

  \indexspace

  \item $Z(A)$, centre of a group $A$, 29

 \end{multicols}\end{theglossary}



\providecommand{\bysame}{\leavevmode\hbox to3em{\hrulefill}\thinspace}
\providecommand{\MR}{\relax\ifhmode\unskip\space\fi MR }
\providecommand{\MRhref}[2]{%
  \href{http://www.ams.org/mathscinet-getitem?mr=#1}{#2}
}
\providecommand{\href}[2]{#2}


\end{document}
